\documentclass{article}
\usepackage{amsfonts}
\usepackage{amssymb}
\usepackage{amsmath}
\usepackage{amsthm}

\newtheorem{theorem}{Theorem}
\newtheorem{proposition}{Proposition}
\newtheorem{example}{Example}
\newtheorem{lemma}{Lemma}
\newtheorem{corollary}{Corollary}

\newcommand{\bbC}{{\mathord{\mathbb{C}}}}
\newcommand{\bbD}{{\mathord{\mathbb{D}}}}
\newcommand{\bbI}{{\mathord{\mathbb{I}}}}
\newcommand{\bbN}{{\mathord{\mathbb{N}}}}
\newcommand{\bbQ}{{\mathord{\mathbb{Q}}}}
\newcommand{\bbR}{{\mathord{\mathbb{R}}}}
\newcommand{\bbT}{{\mathord{\mathbb{T}}}}
\newcommand{\bbZ}{{\mathord{\mathbb{Z}}}}

\newcommand{\ttens}{\mathbin{\widehat\otimes}}
\newcommand{\ev}{\mathop{\mathrm{ev}}\nolimits}
\renewcommand{\Re}{\mathop{\mathrm{Re}}}
\newcommand{\supp}{\mathop{\mathrm{supp}}}
\newcommand{\clos}{\mathop{\mathrm{clos}}}
\newcommand{\const}{\mathop{\mathrm{const}}}
\newcommand{\Diff}{\mathop{\mathrm{Diff}}\nolimits}
\newcommand{\BL}{\mathop{\mathrm{BL}}}
\newcommand{\GL}{\mathop{\mathrm{GL}}}
\newcommand{\Hom}{\mathop{\mathrm{Hom}}\nolimits}

\newcommand{\Sp}{\mathop{\mathrm{Sp}}\nolimits}
\newcommand{\SU}{\mathop{\mathrm{SU}}}
\newcommand{\sll}{\mathop{\mathrm{sl}}}
\newcommand{\SO}{\mathop{\mathrm{SO}}}
\newcommand{\rU}{\mathop{\mathrm{U{}}}}
\newcommand{\Ad}{\mathop{\mathrm{Ad}}}
\newcommand{\Int}{\mathop{\mathrm{Int}}}
\newcommand{\Inn}{\mathop{\mathrm{Inn}}}
\newcommand{\rank}{\mathop{\mathrm{rank}}}

\newcommand{\fin}{{\mathord{\mathrm{fin}}}}
\newcommand{\cM}{\mathord{\mathcal{M}}}
\newcommand{\cB}{\mathord{\mathcal{B}}}
\newcommand{\cI}{\mathord{\mathcal{I}}}
\newcommand{\cJ}{\mathord{\mathcal{J}}}

\newcommand{\cP}{\mathord{\mathcal{P}}}
\newcommand{\cS}{\mathord{\mathcal{S}}}

\newcommand{\frB}{\mathord{\mathfrak{B}}}
\newcommand{\frC}{\mathord{\mathfrak{C}}}
\newcommand{\frF}{\mathord{\mathfrak{F}}}
\newcommand{\frG}{\mathord{\mathfrak{G}}}
\newcommand{\frI}{\mathord{\mathfrak{I}}}
\newcommand{\frH}{\mathord{\mathfrak{H}}}
\newcommand{\frh}{\mathord{\mathfrak{h}}}
\newcommand{\frM}{\mathord{\mathfrak{M}}}
\newcommand{\frP}{\mathord{\mathfrak{P}}}
\newcommand{\frR}{\mathord{\mathfrak{R}}}
\newcommand{\frS}{\mathord{\mathfrak{S}}}
\newcommand{\frg}{\mathord{\mathfrak{g}}}
\newcommand{\frn}{\mathord{\mathfrak{n}}}
\newcommand{\frl}{\mathord{\mathfrak{l}}}
\newcommand{\frL}{\mathord{\mathfrak{L}}}

\newcommand{\frz}{\mathord{\mathfrak{z}}}
\newcommand{\al}{{\mathord{\alpha}}}
\newcommand{\be}{{\mathord{\beta}}}
\newcommand{\ga}{{\mathord{\gamma}}}
\newcommand{\Ga}{{\mathord{\Gamma}}}
\newcommand{\vf}{{\mathord{\varphi}}}
\newcommand{\si}{{\mathord{\sigma}}}
\newcommand{\ka}{{\mathord{\kappa}}}
\newcommand{\vk}{{\mathord{\varkappa}}}
\newcommand{\la}{{\mathord{\lambda}}}
\newcommand{\La}{{\mathord{\Lambda}}}
\newcommand{\de}{{\mathord{\delta}}}
\newcommand{\ep}{{\mathord{\varepsilon}}}
\newcommand{\eps}{{\mathord{\epsilon}}}
\newcommand{\ze}{{\mathord{\zeta}}}
\newcommand{\one}{{\mathord{\mathbf1}}}
\newcommand{\Tr}{\mathop{\mathrm Tr}}
\newcommand{\scal}[2]{\left<#1,#2\right>}
\newcommand{\eqv}[1]{\mathrel{\stackrel{\scriptstyle{#1}}{\sim}}}
\let\wh=\widehat
\let\ov=\overline
\let\td=\tilde
\let\wtd=\widetilde
\newcommand{\<}{\mathrel{\mathchar534}}

\date{}
\title{\bf Maximal ideal spaces of invariant function
algebras on compact groups}
\author{V.M. Gichev}
\begin{document}
\maketitle
\begin{abstract}
The paper contains a description of the maximal ideal spaces
(spectra) $\cM_A$ of bi-invariant function algebras $A$ on a
compact group $G$. There are natural compatible structures in
$\cM_A$: it is a compact topological semigroup with involution,
polar decomposition, and analytic structure. The paper contains a
description of $\cM_A$ and related results on function algebras;
for example, a bi-invariant function algebra on a connected
compact Lie group is antisymmetric if and only if the Haar measure
of its maximal torus is multiplicative on it. Some results are
extended to the case of compact commutative homogeneous spaces. As
a consequence, we get an infinite dimensional version of the
Hilbert--Mumford criterion for commutative (as homogeneous spaces)
orbits of compact connected Lie groups, where one parameter
semigroup is replaced by a finite sequence of them.



\end{abstract}


\section*{Introduction}

Let $G$ be a compact Hausdorff group and $C(G)$ be the Banach
algebra of all continuous functions on $G$ equipped with the
sup-norm. We say that $A$ is an {\it invariant algebra on $G$} if
it is a closed subalgebra of $C(G)$ which is invariant with
respect to all left and right translations. Invariant algebras
could be defined as closed sub-bialgebras of $C(G)$ with identity.
We describe the maximal ideal spaces $\cM_A$ (spectra) of these
algebras. They have a natural semigroup structure (induced by the
convolution of representing measures) and an analytic structure.
The simplest example is the closed unit disc
$\ov\bbD=\bbT\cup\bbD$, where $\bbT$ is the circle group and
$\bbD$ is the open unit disc in the complex plane $\bbC$;
$\ov\bbD$ is the maximal ideal space for the algebra of the
continuous functions on $\bbT$ which admit analytic extensions to
$\bbD$, and a multiplicative subsemigroup of $\bbC$. The
polynomially convex hull $\wh G$ of any compact matrix group $G$
is a semigroup of this type with respect to the matrix
multiplication; for example, $\wh{\rU(n)}$ is the unit matrix ball
in the space of $n\times n$-matrices. If an invariant algebra $A$
is finitely generated (i.e. generated as a Banach algebra by its
finite dimensional bi-invariant subspace), then $\cM_A=\wh G$ for
some compact matrix group $G$. The case of algebras which are not
finitely generated is more complicated, for example, it includes
some algebras of analytic almost periodic functions on the upper
halfplane in $\bbC$.

Roughly speaking, the semigroup $\cM_A$ can be builded by gluing
together complex Lie semigroups whose skeletons are compact
groups. Their units are idempotents in $\cM_A$. There is a natural
involution in $\cM_A$, which extends the inversion in $G$ and
fixes all idempotents. Moreover, there is a polar decomposition:
$\vf=gs$, where $\vf\in\cM_A$, $g\in G$, $s\in\cS_A$, and $\cS_A$
is the set of symmetric nonnegative elements (see (\ref{defsa})
for the precise definition). Each $s\in\cS_A$ is included to a
unique one parameter symmetric semigroup, which extends to the
right halfplane in $\bbC$ and induces the analytic structure in
$\cM_A$.

Let $G$ be abelian and let $S_A$ be the set of those one
dimensional characters of $G$ that are contained in $A$. Then, the
linear span of $S_A$ is dense in $A$. The set $S_A$ is a semigroup
in the dual group and $\cM_A\cong\Hom(S_A,\ov\bbD)$, where $\bbD$
is considered as a multiplicative semigroup with identity and
zero, and the homomorphisms are assumed to be nonzero.

Let $G$ be a Lie group and $T$ be a maximal torus in $G$. Then any
$G$-orbit $\{g^{-1}sg:\,g\in G\}$ of $s\in\cS_A$ intersects the
set $\wh T\mathbin\cap\cS_A$, where $\wh T$ is the $A$-hull of $
T$ (see (\ref{ahull}) for the definition); under some additional
assumptions, this is a maximal abelian subsemigroup of $\cM_A$.

Invariant algebras have been studied since 50s. Initially,
invariant algebras on abelian groups were considered as a natural
generalization of the algebra of analytic and continuous up to the
boundary functions on $\bbD$. In \cite{AS}, Arens and Singer got
some results similar to classical ones, in particular, an analog
of the Poisson integral. They found the realization of $\cM_A$ as
$\Hom(S_A,\ov\bbD)$ and constructed the polar decomposition in
$\Hom(S_A,\ov\bbD)$ (which appeared earlier in Goldman \cite{Go}
and Mackey \cite{Ma}). Paper \cite{LM60} by de Leew and Mirkil
contains basic facts concerning invariant algebras on locally
compact abelian groups (mainly, in the framework of abelian
harmonic analysis). They also proved in \cite{LM61} that any
$\SO(n)$-invariant algebra on a sphere $S^{n-1}$, $n>2$, is
self-adjoint with respect to the complex conjugation. Similar
results were obtained by Wolf \cite{Wo} and Gangolli \cite{Ga}.
Wolf characterized compact groups $G$ having the property that
each invariant algebra on $G$ is self-adjoint; Gangolli
independently proved that this property holds for connected
compact semisimple Lie groups. Due to the Stone--Weierstrass
theorem, if $A$ is self-adjoint, then it consists of all
continuous functions on $G$ which are constant on cosets of some
normal subgroup of $G$. Here is Wolf's condition: the image of any
one dimensional character of $G$ is finite. As an easy consequence
of the results of this paper, we get another version of the
criterion: the centre of $G$ is profinite (hence, the two
conditions are equivalent, but it is not difficult to check the
equivalence directly). Further, Rider \cite{Ri} proved that a
compact group which admits an antisymmetric Dirichlet invariant
algebra is connected and abelian (by definition, a Dirichlet
function algebra has no nontrivial real orthogonal measure; an
antisymmetric function algebra contains no nonconstant real
functions). Antisymmetric invariant algebras were characterized by
Rosenberg in \cite{Ros} by three conditions which are too
complicated to be formulated here; \cite{Ros} also contains a
generalization of \cite{Ri}. It was noted in \cite{Gl} and
\cite{Gi79} that an invariant algebra $A$ is antisymmetric if and
only if the Haar measure is multiplicative on $A$ (see Section~4).
Also, the paper \cite{Gl} contains a construction for analytic
discs in the maximal ideal spaces of invariant algebras on locally
compact abelian groups.

The approach based on the observation that $\cM_A$ has a natural
semigroup structure was used in \cite{Gi79} where the case of
noncompact groups was mainly considered (\cite{Gi79} contains an
error concerning  the structure of Lie algebras admitting an
invariant cone which is essential only for noncompact groups and
does not affect the main results). The results on compact groups
were announced in \cite{Gi79faa}; their proofs were published in
\cite{Gi80}. This paper, in particular, contains a new exposition
of an essential part of these results (\cite{Gi80} is hardly
accessible).

There are several fields which are closely connected with the
object of this paper. An evident one is the class of invariant
algebras on homogeneous spaces (in particular, algebras of
CR-functions; this involves the problem of description of the
polynomially convex hulls for orbits of compact groups, etc.).
This case essentially differs from the case of bi-invariant
algebras (see \cite{La}, \cite{GL}). The last section contains a
partial extension of the results to invariant algebras on
commutative (spherical, multiplicity free) homogeneous spaces (the
group $G$ is spherical as the homogeneous space $G\times G/G$). As
a consequence, we get a version of the Hilbert--Mumford criterion
for Banach spaces and spherical orbits.

The paper is organized as follows.

Section~1 contains the notation, the definitions and some known
facts.

In Section~2, the basic properties of $\cM_A$ as a semigroup with
involution are found. It can be realized as the semigroup of those
endomorphisms of $A$ which commute with all left translations.
These automorphisms are bounded with respect to $L^2(G)$-norm.
Hence they extend to the closure $H^2$ of $A$ in $L^2(G)$, and the
involution agrees with the conjugation for operators in the
Hilbert space $H^2$. For abelian $G$, there is another natural way
to define the structure of a semigroup with an involution in
$\cM_A$; it is proved that the two structures coincide. Also, we
show that an invariant algebra is the same as a closed
sub-bialgebra with identity of the Hopf algebra $C(G)$.

The closure of the restriction of an invariant algebra to a closed
subgroup is an invariant algebra. In Section~3, its maximal ideal
space is identified with the $A$-hull of this subgroup. The
averaging by left over a closed subgroup $H\subseteq G$ defines a
projection to the algebra $A^H$ of all left $H$-invariant
functions in $A$. It is proved that the dual map
$\cM_A\to\cM_{A^H}$ is surjective. Moreover, the semigroup $\cM_A$
acts by right on $\cM_{A^H}$ in such a way that the orbit of any
point in the Shilov boundary coincides with $\cM_{A^H}$. The
averaging over normal subgroups makes possible to realize the
compact semigroup $\cM_A$ for an invariant algebra $A$ on an
arbitrary compact group $G$ as the inverse limit of maximal ideal
spaces of invariant algebras on Lie groups.

In Section~4, the results of the previous section are applied to
the studying of idempotents in $\cM_A$. It is proved that $A$ is
antisymmetric if and only if the Haar measure on $G$ is
multiplicative on $A$, and this is equivalent to the condition
that $\cM_A$ is a semigroup with zero. The set of idempotents is a
complete lattice with respect to the natural order in it.

The polar decomposition in $\cM_A$ and one dimensional analytic
structure are constructed in Section~5.

Section~6 contains the main results of the paper. We assign  to
each idempotent $j$ several objects: a group $G^j\subseteq\cM_A$
with the identity element $j$, a cone $C^j$ in its Lie algebra
$\frg^j$, sets $S^j=\exp(iC^j)$ and $P^j=G^jS^j$, etc. The cone
$C^j$ is pointed, closed, convex, $\Ad(G^j)$-invariant, the set
$P^j$ is a semigroup in the complexification of $G^j$, and $P^j$
is naturally embedded to $\cM_A$. This defines the analytic
structure in $\cM_A$.  Any idempotent in $\cM_A$ is joined with
the unit of $G$ by a chain of one parameter semigroups. Also, the
section contains a criterion of the antisymmetry for the invariant
algebras on connected Lie groups, which is stronger than the
criterion of Section~4: $A$ is antisymmetric if and only if the
Haar measure of a maximal torus $T$ is multiplicative on $A$
(equivalently, if the restriction of $A$ to $T$ is antisymmetric).
The $A$-hull $\wh T$ of $T$ is an abelian subsemigroup of $\cM_A$
such that each symmetric $s\in\cM_A$ is conjugated to some element
of $\wh T$; moreover, under some additional assumptions, $\wh T$
is a maximal abelian subsemigroup of $\cM_A$.

In Section~7, the case of tori is considered. Then
$\cM_A=\Hom(S,\ov\bbD)$, where $S=A\mathbin\cap \wh G$ and $\wh G$
is the dual to $G$ group, which is embedded to $C(G)$. The major
point is the identification of idempotents. Any idempotent
corresponds to the characteristic function of a semigroup
$P\subseteq S$ such that $S\setminus P$ is a semigroup ideal. They
can be found by a procedure, which is described in this section.
If $A$ is finitely generated, then the procedure consists of only
one step.

Section~8 contains several illustrating examples and some
additional results. In particular, it is proved that each
invariant algebra on $G$ is self-adjoint if and only if the centre
of $G$ is profinite; also, the section contains a description of
antisymmetric Dirichlet invariant algebras on Lie groups.

In Section~9, we consider invariant algebras on spherical
homogeneous spaces of connected compact Lie groups (they are also
known as commutative or multiplicity free homogeneous spaces).
This class includes the groups $G$ considered as homogeneous
spaces of $G\times G$, where $G$ acts on itself by the left and
right translations. They has been studied last decade from various
viewpoints. Some of the results above hold for invariant algebras
on spherical homogeneous spaces (for instance, an invariant
algebra is antisymmetric if and only if the normalized invariant
measure is multiplicative; for a generic homogeneous space, this
is not true).
Using them, we get an infinite dimensional version of the
Hilbert--Mumford criterion for the spherical orbits of compact
groups. In it, one parameter semigroup is replaced by a finite
sequence of them. An example shows that the finite dimensional
version does not hold for Banach spaces, in any sense.


\section{Preliminaries}
Everywhere in this paper, $G$ denotes a compact Hausdorff
topological group with identity element $e$, and $C(G)$ is the
space of all continuous functions on $G$ endowed with the
$\sup$-norm, which is denoted by $\|\kern5pt\|$ (other norms are
indicated explicitly). For function spaces $L$ and $M$, $LM$ or
$L\cdot M$ is the linear span of products $uv$, where $u\in L$,
$v\in M$. Similarly, if $A,B\subseteq G$, then $AB=A\cdot
B=\{ab:\,a\in A,\, b\in B\}$. Let $\widehat G$ be the dual object
to $G$, i.e., the set of classes of equivalent irreducible complex
linear representations of $G$. For a finite dimensional
representation $\tau$, $M_\tau$, $V_\tau$, $\scal{\ }{\ }_\tau$,
$\chi_\tau(g)=\Tr\tau(g)$ denote the space of the matrix elements,
the space of the representation, the invariant inner product in
it, and the character, respectively. The following equalities hold
for all finite dimensional representations $\tau,\sigma$:
\begin{equation}\label{mtms}
M_{\tau\otimes\sigma}=M_\tau M_\sigma,\quad
\chi_{\tau\otimes\sigma}=\chi_\tau\chi_\sigma.
\end{equation}
If $\tau\in\wh G$ and $\si$ is a continuous representation of $G$
(in general, infinite dimensional), then
$P_\tau=\int_G\si(g)\chi_\tau(g)\,dg$, where $dg$ denotes the Haar
measure in $G$, is the projection to the $\tau$-isotypical
component in $V_\si$. The convolution $\mu*\nu$ of finite regular
Borel measures $\mu,\nu$ on $G$ can be defined by the equality
\begin{eqnarray}\label{mudef}
&\int_G f(g)\,d\mu*\nu(g)=\int_{G\times G}f(gh)\,d\mu(g)d\nu(h),
\end{eqnarray}
where the space $M(G)$ of these measures is identified with the
dual to $C(G)$ Banach space. Equipped with ${}*{}$, $M(G)$ is a
Banach algebra. The spaces $C(G)$, $L^p(G,dg)$, $p\geq1$,
naturally embedded to $L^1(G)\subseteq M(G)$, are subalgebras of
$M(G)$ (in the sequel, we denote $L^p(G,dg)$ by $L^p(G)$). Spaces
$M_\tau$, $\tau\in\wh G$, are their minimal ideals ("ideal" means
"two-sided ideal"). In particular,
\begin{eqnarray}\label{cotri}
&M_\tau*M_\si=\left\{\begin{array}{cc}M_\tau,&\tau=\si,\\
0,&\tau\neq\si.\end{array}\right.
\end{eqnarray}
By Peter--Weyl Theorem,
\begin{eqnarray}\label{l2dec}
&L^2(G)=\sum_{\tau\in\wh G}\oplus M_\tau,
\end{eqnarray}
where the sum is orthogonal with respect to the standard inner product
\begin{eqnarray}\label{scal2}
\scal{u}{v}=\int_Gu(g)\ov{v(g)}\,dg
\end{eqnarray}
and the bar denotes the complex conjugation. Throughout the paper
(except for the last section), a function space on $G$ is called
{\it invariant} if it is bi-invariant, i.e. invariant with respect
to left and right translations defined by
\begin{eqnarray*}
L_gf(h)=f(gh),\qquad R_gf(h)=f(hg),
\end{eqnarray*}
respectively. Note that $L$ is not a representation. For a
function space $F$, set
\begin{eqnarray}
&\Sp F=\{\tau\in \widehat G:\,M_\tau\subseteq F\},\nonumber \\
&F_\fin=\sum\nolimits_{\tau\in\Sp F} M_\tau,\label{lfind}
\end{eqnarray}
where the sum is algebraic (i.e., the set of finite sums of
vectors in summands).

Throughout the paper, ``weak" refers to the $*$-weak topology of
the dual space, and ``strong" to the strong operator topology in
the space $\BL(X)$ of bounded linear operators in a Banach space
$X$. The dual to $X$ space is denoted by $X'$. All representations
are supposed to be strongly continuous. For a representation
$\rho$, we denote by $\Sp\rho$ the set of its irreducible
components (without multiplicities). The proof of the following
proposition is omitted since these facts are well-known (see, for
example, \cite{HR}).
\begin{proposition}\label{inveq}
A closed subspace $F$ of $C(G)$, or $L^p(G)$, $1\leq p<\infty$, or
weakly closed subspace of $M(G)$, is an ideal of the convolution
algebra $M(G)$ if and only if it is invariant. Furthermore, the
space $F_\fin$ is dense in $F$ for any closed invariant space $F$
(in the norm topology if $F\subseteq L^p(G)$ or $F\subseteq C(G)$
and weakly if $F\subseteq M(G)$) and
\begin{eqnarray}\label{lfinc}
F_\fin=F\cap C(G)_\fin.
\end{eqnarray}
\end{proposition}
\noindent A simple computation with matrix elements shows that
each space $M_\tau$, $\tau\in\wh G$, is invariant with respect to
the antilinear involution
\begin{eqnarray}\label{infde}
f^\star(g)=\ov{f(g^{-1})}.
\end{eqnarray}
Proposition~\ref{inveq} and {\rm(\ref{lfind})} imply the following
corollary.
\begin{corollary}\label{l*=l}
$F^\star=F$ for any $F$ as above.\qed
\end{corollary}
The space $F_\fin$ could be defined as the set of all functions in
$F$ such that the linear span of their left and right translates
is finite-dimensional.

The {\it maximal ideal space} $\cM_A$ ({\it spectrum}) of a
commutative Banach algebra $A$ with identity element $\one$ is the
set of all nonzero multiplicative linear functionals (i.e.
$\cM_A=\mathop{\rm Hom}(A,\bbC)$) equipped with the weak topology
of the dual space $A'$.
For a Hausdorff compact space $Q$, $M(Q)$ denotes the dual to
$C(Q)$ space of all finite regular complex measures on $Q$. If
$\vf\in\cM_A$, then $\|\vf\|_{A'}=1$. Since $\vf(\one)=1$, any
norm preserving extension of $\vf$ from $A$ to $C(Q)$ is a
positive measure of the total mass 1. These measures are called
{\it representing}; the set of all representing measures for
$\vf\in\cM_A$ is denoted by $\cM_\vf$. This is a weakly compact
convex subset of $M(Q)$. A measure that is representing for some
$\vf\in\cM_A$ will be called {\it multiplicative}. Every $x\in Q$
corresponds the evaluation functional
\begin{eqnarray*}
\ev_x(f)=f(x).
\end{eqnarray*}
The Dirac measure at $x$ is representing for $\ev_x$. The image of
$Q$ under the mapping $\ev:\,x\to \ev_x$ can be identified with
the factor of $Q$ by the following equivalence:
\begin{eqnarray*}
x\eqv{A}y\quad\Longleftrightarrow\quad\ f(x)=f(y)\quad
\mbox{for all}\quad f\in A.
\end{eqnarray*}
If each class is a single point, then $A$ is said to be {\it
separating}.  In this case, $\ev$ is an embedding and $Q$ may be
considered as a subset of $\cM_A$. Let $A$ be a closed separating
subalgebra of $C(G)$. The Gelfand transform $A\to C(\cM_A)$ is
defined by $\widehat f(\vf)=\vf(f)$. In most cases, we omit the
hat.

We say that $A\subseteq C(G)$ is an {\it invariant algebra}, if it
is a closed invariant subalgebra of $C(G)$ that contains constant
functions. If $A$ is an invariant algebra on the group $G$, then
the relation $\eqv A$ above is invariant. Hence, the equivalence
class of the identity $e$ is a closed normal subgroup $N$ and the
algebra $A$ may be considered as a separating invariant  algebra
on the group $G/N$.
Clearly, an invariant subspace $A\subseteq C(G)$ is an algebra if
and only if $A_\fin$ is an algebra. It follows from
{\rm(\ref{mtms})} that $A_\fin$ is an algebra if and only if $\Sp
A$ contains irreducible components of $\tau\otimes\sigma$ for
every $\tau,\sigma\in\Sp A$.

Let $S$ be a semigroup or an algebra with unit $e$. Then, $x\in S$
is said to be {\it invertible}, if there exist $y,z\in S$ such
that $xy=zx=e$ (then $y=z$). The set of all invertible elements in
$S$ is a group, which we denote by $S^{-1}$.

Throughout the paper (except for the last section), $A$ is an
invariant algebra on $G$, and $H^2$ is its closure in $L^2(G)$. It
follows from {\rm(\ref{cotri})}, {\rm(\ref{l2dec})}, and
Proposition~\ref{inveq} that
\begin{eqnarray*}
H^2=\sum\nolimits_{\tau\in\Sp A}\oplus M_\tau,
\end{eqnarray*}
where the sum is orthogonal.

Let $\mu\in M(Q)$, where $Q$ is a compact Hausdorff space, $\phi$,
$\Phi$ be a continuous mapping from $Q$ to a Banach space $X$ and
a strongly continuous mapping $\Phi:\,Q\to\BL(X)$, respectively.
Then their integrals are defined by
\begin{eqnarray}
\xi\Big(\int_Q\phi(q)\,d\mu(q)\Big)=
\int_Q\xi\big(\phi(q)\big)\,d\mu(q) \quad\text{for all}\ \xi\in X',\\
\label{intd1}
\Big(\int_Q\Phi(q)\,d\mu(q)\Big)x=
\int_Q\Phi(q)x\,d\mu(q) \quad\text{for all}\ x\in X.
\label{intd2}
\end{eqnarray}
The latter formula extends a strongly continuous representation $\tau$
to the convolution algebra $M(G)$:
\begin{eqnarray}\label{exten}
\tau(\mu)=\int_G\tau(g)\,d\mu(g).
\end{eqnarray}

The unit circle and the open (closed) unit disc in $\bbC$ are
denoted by $\bbT$, $\bbD$ ($\ov\bbD$), respectively. In most
cases, they are considered as subsemigroups of the multiplicative
semigroup $\bbC$. Also, $\bbC^+$ (the closed right halfplane in
$\bbC$) and $\bbR^+=[0,\infty)$ are equipped with the structure of
additive semigroup.

We need some basic facts on function algebras; for more details,
see \cite[Ch. 2]{Gam}. Let $Q$ be compact and Hausdorff,
$A\subseteq C(Q)$ be a closed subalgebra. The set $E\subseteq Q$
is called a {\it peak set for $A$} if there exists $f\in A$ such
that $f(x)=1$, $x\in E$, and $|f(x)|<1$, $x\notin E$. We say that
$f$ is a {\it peak function for $E$}. A {\it peak point} is a peak
set consisting of a single point. A {\it $p$-set} is a set which
can be realized as the intersection of a family of peak sets; a
{\it  $p$-point} is a point with this property. Clearly, $p$-sets
are closed. If $A$ is separating, then for each $f\in A$ the
function $|f|$ attains its maximal value at some $p$-point. In
particular, the set of these points is not empty. The Dirac
measure at a $p$-point $x$ is the unique representing measure for
the $\ev_x$. A metrizable $Q$ contains a peak point. Modulo these
facts, the following lemma is obvious.
\begin{lemma}\label{peakr}
Let $A$ be a separating $G$-invariant algebra on a homogeneous
space $M=G/H$ of a compact group $G$.  Then each point of $M$ is a
$p$-point. Moreover, if $G$ is a Lie group then all points of $M$
are peak points.
\end{lemma}
Let $E\subseteq Q$ be closed and $B$ be the closure in $C(E)$ of
the restriction $A|_E$. Then every $\vf\in\cM_B$ can be considered
as a functional on $A$. This defines the embedding
$\cM_B\to\cM_A$. Its image consists of $\vf\in\cM_A$ that admit
continuous extension to $B$. This set is called the {\it $A$-hull}
of $E$. It is denoted by $\wh E$, and can also be defined by
\begin{eqnarray}
\label{ahull}
\wh E=\{\vf\in\cM_A:\ |\vf(f)|\leq\sup_{x\in E}|f(x)|
\quad\text{for all}\ f\in A\}.
\end{eqnarray}
Since $A$ can be considered as a subalgebra of $C(\cM_A)$, the
definition can be applied to subsets of $\cM_A$. A set
$E\subseteq\cM_A$ is called {\it $A$-convex} if $\wh E=E$.
Replacing $A$ by the algebra $\cP(V)$ of all holomorphic
polynomials on a complex linear space $V$, and $\cM_A$ by $V$, we
get the definition of the {\it polynomially convex hull} of a
compact set $E\subset V$.
\begin{lemma}\label{rpeak}
Let $E$ be a closed subset of $Q$. Then $\vf\in\wh E$ if and only
if there exists $\mu\in\cM_\vf$ such that $\supp\mu\subseteq E$.
If $E$ is a $p$-set and $\vf\in\wh E$, then $\supp\mu\subseteq E$
for any $\mu\in\cM_\vf$.
\end{lemma}
\begin{proof}
If $\mu\in\cM_\vf$ and $\supp\mu\subseteq E$, then the inequality
$|\vf(f)|\leq\|f\|$ is obvious. Hence, $\vf\in\wh E$. If
$\vf\in\wh E$, then {\rm(\ref{ahull})} and the Hahn--Banach
theorem imply the existence of a norm preserving extension $\mu$
of $\vf$ to $C(E)$. Since $\mu(E)=\vf(1)=1$, this proves the first
assertion. Proving the second, we may assume without loss of
generality that $E$ is a peak set. Let $f$ be a peak function for
$E$. Then $\vf(f)=1$, since $f|_E$  is identity element of the
restricted algebra. If $x\notin E$, then $|f(x)|<1$ by definition,
hence, the assumption $\supp\mu\not\subseteq E$ implies $|\int
f\,d\mu|<1$.
\end{proof}

We do not use any deep theorem of the theory of topological
semigroups or measure convolution semigroups; it is sufficient to
mention following elementary facts: 1) each compact topological
semigroup contains the unique minimal ideal which is the union of
groups; 2) each convolution-idempotent positive measure on $G$ of
total mass 1 is the Haar measure of some closed subgroup. Proofs
can be found in \cite{HM} and \cite{He}, respectively. An element
$\eps$ of a semigroup $S$ is called {\it zero} if $\eps
s=s\eps=\eps$ for all $s\in S$. For a subset $X\subseteq S$ let
\begin{eqnarray*}
&Z_S(X)=\{z\in S:\,zx=xz\ \text{for all}\ x\in X\},\\
&N_S(X)=\{z\in S:\,zX=Xz\}.
\end{eqnarray*}
If $S=\cM_A$, then the index will be omitted. One parameter
semigroup in a topological semigroup $S$ is a continuous
homomorphism $\bbR^+\to S$.

The {\it analytic structure} in $\cM_A$ is a mapping
$\la:\,D\to\cM_A$, where $D$ is a domain in some complex manifold,
such that $f(\la(z))$ is analytic for all $f\in A$.

Lie algebras of Lie groups are realized as algebras of left
invariant vector fields and are denoted by the corresponding
lowercase Gothic letters.

The {\it relative interior} of a set $X$ in a real vector space $V$
is its interior in its linear span. It will be denoted by $\Int X$.
A cone $C\subseteq V$ is a subsemigroup of the additive group of $V$
which is invariant with respect to dilations $v\to \la v$, $\la>0$
(hence we consider only convex cones). There is the natural preorder
in a cone $C$:
\begin{eqnarray}\label{preco}
u\< v\ \Longleftrightarrow\ -\ep u+(1+\ep)v\in C\ \mbox{for some}\  \ep>0.
\end{eqnarray}
It defines the equivalence whose
classes are called {\it faces}.
For $x\in C$, we denote by $F_x$ or $F_{x,C}$ the face that contains $x$.
Clearly, each face is open in its linear span, the closure of $F_x$ in $C$
coincides with the union of those faces of $C$ which are contained
in the linear span of $F_x$, and with the set
\begin{eqnarray}
\ov F_x=\{y\in C:\,y\< x\}. \label{clofa}
\end{eqnarray}
We shall say that $\ov F_x$ is a {\it closed face}
and denote by $\frF_C$ the family of all closed faces of $C$.

A closed cone $C$ is called {\it pointed} if $C\mathbin\cap(-C)=0$
(in other words, if $\{0\}$ is a face of $C$). An arbitrary cone
is pointed, if its closure has this property.

Condition {\rm(\ref{preco})} for cones is evidently equivalent to
the following one:
\begin{eqnarray}\label{prela}
u\< v\ \Longleftrightarrow\ -ku+lv\in C\ \mbox{for some}\
k,l\in\bbN,
\end{eqnarray}
where $\bbN= \{1,2,\dots\}$. The definition {\rm(\ref{prela})},
hence {\rm(\ref{clofa})}, can be applied to semigroups in
$\bbZ^n$.
The group $\bbZ^n$ is always assumed to be canonically embedded to
$\bbR^n$. Any semigroup $S\subseteq\bbZ^n$ has the {\it asymptotic
cone}
\begin{eqnarray}\label{acone}
&\al(S)=\clos\bigcup\limits_{n=1}^\infty\frac1n S\subseteq\bbR^n.
\end{eqnarray}
Clearly, $\al(S)$ is a closed convex cone, and $S\subseteq\al(S)$.
Being convex, $\al(S)$ coincides with its bi-dual; the dual to a
set $X$ cone $X^\star$ is defined by
\begin{eqnarray*}
X^\star=\{y\in\bbR^n:\,\scal{x}{y}\geq0\ \text{for all}\ x\in X\},
\end{eqnarray*}
where $\scal{\ }{\ }$ is the standard inner product in $\bbR^n$.
Obviously, $X^\star=\left(\clos X\right)^\star$ and
$X^\star=\left(\frac1n X\right)^\star$ for all $n\in\bbN$. Therefore,
$S^\star=\al(S)^\star$, and this gives another definition of $\al(S)$:
\begin{eqnarray*}
\al(S)=S^\star{}^\star.
\end{eqnarray*}

\section{Semigroup structure in maximal ideal spaces}
The algebraic tensor product $C(X)\otimes C(Y)$ (over $\bbC$),
where $X,Y$ are Hausdorff compact spaces, can be identified with
the linear span of functions $f(x)g(y)$ on $X\times Y$, where
$f\in C(X)$, $g\in C(Y)$. Thus the space $C(X\times Y)$ can be
considered as a completion of $C(X)\otimes C(Y)$. We shall denote
it by $C(X)\ttens C(Y)$ and identify with $C(X\times Y)$. If
$X=Y=G$, where $G$ is a compact group, then the group
multiplication $G\times G\to G$ induces the {\it comultiplication}
$\iota:\,C(G)\to C(G)\ttens C(G)=C(G\times G)$ by
\begin{eqnarray}
\label{iodef}
\iota(f)(g,h)=f(gh),\quad f\in C(G).
\end{eqnarray}
The space $C(G)$ with the pointwise multiplication of functions,
this comultiplication, unit $\one$, counit $\de_e$, and the
antipode mapping $f(g)\to f(g^{-1})$ is a Hopf algebra. In any
bialgebra $X$, the formula
\begin{eqnarray}\label{convo}
&\mu*\nu(f)=(\mu\otimes\nu)(\iota(f)),\quad\mu,\nu\in X',\ f\in X
\end{eqnarray}
defines a multiplication in $X'$. If $X=C(G)$, then
{\rm(\ref{convo})} coincides with {\rm(\ref{mudef})}. For a
subspace $L\subseteq C(G)$, let $L\ttens L$ be the closure of
$L\otimes L$ in $C(G)\ttens C(G)$. We say that $A\subseteq C(G)$
is a {\it sub-bialgebra} if
\begin{eqnarray}
&\one\in A,\label{coone}\\
&A\cdot A\subseteq A,\label{hdefa}  \\
&\iota(A)\subseteq A\ttens A
\label{hdefh}.
\end{eqnarray}
Clearly, $\iota$ is isometric. By {\rm(\ref{iodef})}
and {\rm(\ref{convo})},
\begin{eqnarray}\label{baain}
&\|\psi*\vf\|\leq\|\psi\|\|\vf\|.
\end{eqnarray}
\begin{proposition}\label{hosub}
A closed subspace $A\subseteq C(G)$ is an invariant algebra if and
only if it is a sub-bialgebra of $C(G)$. Then $A^\bot$ is an ideal
in $M(G)$ and the algebra $A'$ with the multiplication
{\rm(\ref{convo})} is isomorphic to $M(G)/A^\bot$.
\end{proposition}
\begin{proof}
Since {\rm(\ref{coone})} and {\rm(\ref{hdefa})} is exactly the
assumption that $A$ is an algebra with the unit $\one$, it is
sufficient to prove that {\rm(\ref{hdefh})} holds if and only if
$A$ is invariant. It follows from {\rm(\ref{mudef})} and
{\rm(\ref{scal2})} that $\scal{u\otimes
v}{\iota(w)}=\scal{u*v}{w}$ for all $u,v,w\in C(G)_\fin$. By
{\rm(\ref{cotri})}, {\rm(\ref{l2dec})} and {\rm(\ref{convo})},
\begin{eqnarray}\label{mtm}
&\iota(M_\tau)\subseteq M_\tau\otimes M_\tau,\quad \tau\in\wh G.
\end{eqnarray}
If $A$ is invariant, then {\rm(\ref{mtm})}, taken together with
Proposition~\ref{inveq} (particularly {\rm(\ref{lfinc})}), implies
{\rm(\ref{hdefh})}. Conversely, {\rm(\ref{hdefh})},
{\rm(\ref{convo})} and obvious relations
\begin{eqnarray*}
&A^\bot\otimes M(G)\perp A\otimes C(G),\quad
M(G)\otimes A^\bot\perp C(G)\otimes A
\end{eqnarray*}
imply that $A^\bot$ is an ideal in $M(G)$. Since $A^\bot$ is
weakly closed, it is invariant by Proposition~\ref{inveq}. Hence,
the same is true for $A$. It follows from {\rm(\ref{mudef})},
{\rm(\ref{convo})}, {\rm(\ref{hdefh})} that the restriction of
linear functionals to $A$ is a homomorphism $M(G)\to A'$.
\end{proof}
By Corollary~\ref{l*=l}, $A^\star=A$;  evidently, ${}^\star$
is an antilinear automorphism of $A$. It induces an involution
in $\cM_A$ by
\begin{eqnarray}\label{invma}
\vf^*(f)=\ov{\vf(f^\star)},
\end{eqnarray}
where the bar denotes the complex conjugation.
\begin{theorem}\label{semit}
Let $A$ be a separating invariant algebra on a compact group $G$.
Then the space $\cM_A$ with the multiplication {\rm(\ref{convo})}
is a compact topological semigroup. The mapping $\ev$ realizes an
isomorphic embedding of $G$ to $\cM_A$. The identity element $e\in
G$ is the unit of $\cM_A$, and $\cM_A^{-1}=G$; moreover,
conditions $\vf,\psi\in\cM_A$ and $\vf*\psi\in G$ imply
$\vf,\psi\in G$. The involution {\rm(\ref{invma})} is an
involutive antiautomorphism of the semigroup $\cM_A$, which
coincides with the inversion $g\to g^{-1}$ on $G$.
\end{theorem}
\begin{lemma}\label{homeo}
For each $\vf\in A'$ there exists the unique continuous linear
operator $R_\vf:\,A\to A$ satisfying the equality
\begin{eqnarray}\label{defrf}
\psi(R_\vf f)=\psi*\vf(f)
\end{eqnarray}
for all $\psi\in A'$ and $f\in A$. The mapping $R:\,\vf\to R_\vf$
is the extension of the right regular representation $R$ which can
also be defined by {\rm(\ref{exten})}; precisely,
\begin{eqnarray}\label{intrf}
R_\vf=\int_GR_g\,d\mu(g)
\end{eqnarray}
for any extension $\mu\in M(G)$ of the functional $\vf$ from $A$ to $C(G)$.
Furthermore,
$\|R_\vf\|\leq\|\vf\|$,
$R_\vf$ commutes with all $L_g$, $g\in G$,
and the mapping $R:\,\vf\to R_\vf$
restricted to the unit ball $\cB\subset A'$ is a homeomorphism between
$\cB$ with the weak topology and the set
\begin{eqnarray*}
\frB=\{T\in\BL(A):\,\|T\|\leq1,\  TL_g=L_gT\ \mbox{for all}\ g\in
G\}
\end{eqnarray*}
with the strong operator topology.
The inverse mapping $R^{-1}$ corresponds to each $T\in\frB$ the functional
\begin{eqnarray}\label{inver}
R^{-1}(T):\,f\to Tf(e).
\end{eqnarray}
\end{lemma}
\begin{proof}
For any fixed $\vf$ and $f$, the right side of {\rm(\ref{defrf})}
is weakly continuous on $\psi$ due to {\rm(\ref{convo})}: clearly,
$\psi\otimes\vf(u\otimes v)$ is weakly continuous on $\psi$,
hence, this is true for finite sums $\sum u_k\otimes v_k$ and
their limits in the norm topology, in particular, for $\iota(f)$.
Thus, {\rm(\ref{defrf})} correctly defines $R_\vf f\in A$.
Obviously, $R_\vf$ is linear. Since $A^\bot$ is an ideal, it
follows from {\rm(\ref{defrf})}, {\rm(\ref{intd1})} and
{\rm(\ref{intd2})} that the integral in {\rm(\ref{intrf})} is
independent of the choice of $\mu$, and the simple calculation
with {\rm(\ref{convo})} shows that {\rm(\ref{defrf})} holds for
the right side of {\rm(\ref{intrf})}. The inequality
{\rm(\ref{baain})} implies $\|R_\vf\|\leq\|\vf\|$. Using the
approximation of $\iota(f)$ as above and applying
{\rm(\ref{mudef})}, {\rm(\ref{iodef})} and
Proposition~\ref{hosub}, for all $g,h\in G$, $\vf\in A'$ and $f\in
A$, we get
$\de_g*\vf(f)=\de_g\otimes\vf(\iota(f))=\vf(L_gf)$. 
By {\rm(\ref{defrf})},
\begin{eqnarray}\label{rfvfg}
R_\vf f(g)=\de_g(R_\vf f)=\de_g*\vf(f)=\vf(L_gf).
\end{eqnarray}
Hence, $R_\vf L_hf(g)=\vf(L_gL_hf)=\vf(L_{hg}f)=R_\vf
f(hg)=L_hR_\vf f(g)$, and $R_\vf$ commutes with all $L_h$, $h\in
G$. Pick $\ep>0$ and $f\in A$. There exists a neighbourhood $U\ni
e$ and $g_1,\dots,g_n\in G$ such that $\|L_hf-f\|<\ep$ for all
$h\in U$ and
\begin{eqnarray}\label{tmpun}
G=\bigcup_{k=1}^nUg_k.
\end{eqnarray}
Inequalities $|(\psi-\vf)(L_{g_k}f)|<\ep$, $k=1,\dots,n$, define a
weak neighbourhood $V$ of $\vf\in\cB$. If $g\in Ug_k$ then
$g=hg_k$ for some $h\in U$ and for each $\psi\in V\mathbin\cap\cB$
\begin{eqnarray*}
|R_\psi f(g)-R_\vf f(g)|=|(\psi-\vf)(L_{g}f)|\leq
|\psi((L_{hg_k}-L_{g_k})f)|+\\
|(\psi-\vf)(L_{g_k}f)|+
|\vf((L_{hg_k}-L_{g_k})f)|<
|\psi(L_{g_k}(L_{h}f-f))|+\ep \\ +|\vf(L_{g_k}(L_{h}f-f))|
\leq\|L_{h}f-f\|+\ep+\|L_{h}f-f\|<3\ep.
\end{eqnarray*}
By {\rm(\ref{tmpun})}, this is true for all $g\in G$, thus
$\|R_\vf f- R_\psi f\|\leq3\ep$. Therefore, $R$ is continuous.
Due to the inequality $\|R_\vf\|\leq\|\vf\|$, $R$ maps $\cB$ to $\frB$.
Since $\cB$ is compact,
$R$ is a homeomorphism if the inverse mapping exists; to prove this,
note that $R_\vf(e)=\vf(f)$ by {\rm(\ref{defrf})} and {\rm(\ref{rfvfg})},
and that the norm of the functional $f\to Tf(e)$ is not greater than $\|T\|$.
\end{proof}
Let $\frM_A$ be the set of all continuous nonzero endomorphisms of
the algebra $A$ that commute with left translations. We endow
$\frM_A$ with the strong operator topology. Clearly, $\frM_A$ is a
semigroup and $\frM_A\subset\frB$.
A representation $\tau$ of the semigroup $\cM_A$ in a Hilbert
space $H$ is said to be a ${}^*$-representation if
$\tau(\vf^*)=\tau(\vf)^*$ for all $\vf\in\cM_A$.
\begin{theorem}\label{mabal}
The mapping $R$ is a topological isomorphism between $\cM_A$ and $\frM_A$.
For each $\vf\in\cM_A$, $R_\vf$ admits the unique continuous extension
to $H^2$. This defines a ${}^*$-representation $R$ of $\cM_A$ that
realizes a homeomorphic and isomorphic embedding of $\cM_A$
to $\BL(H^2)$.
\end{theorem}
\begin{proof} Due to Lemma~\ref{homeo}, it is sufficient to
prove that $R(\cM_A)=\frM_A$ to check the first assertion of the
theorem but this is an easy consequence of {\rm(\ref{rfvfg})} and
{\rm(\ref{inver})}. Formula {\rm(\ref{intrf})} defines an
extension of $R$ to $H^2$, which is evidently strongly continuous
and unique. Since $R$ is one-to-one on $\cM_A$ and $\cM_A$ is
compact, $R$ is homeomorphic. It follows from {\rm(\ref{infde})}
and {\rm(\ref{invma})} that $\de_g^*=\de_{g^{-1}}$ for all $g\in
G$; hence, for $\vf=\ev_g$ we have
\begin{eqnarray}\label{conjs}
R_{\vf^*}=R_\vf^*,
\end{eqnarray}
where $R_\vf^*$ is the adjoint of $R_\vf$ in the Hilbert space
$H^2$. Since any representing measure is positive and the mappings
$\vf\to\vf^*$ and $R_\vf\to R_\vf^*$ are antilinear and
continuous, (\ref{intrf}) implies {\rm(\ref{conjs})} for
$\vf\in\cM_A$.
\end{proof}
\begin{corollary}\label{reext}
Let $\rho$ be a strongly continuous representation of the group
$G$ in a Banach space. If $\Sp\rho\subseteq\Sp A$, then formula
{\rm(\ref{exten})}, with $\mu\in\cM_\vf$, defines an extension of
$\rho$ to a strongly continuous representation of $\cM_A$. If
$\rho$ is unitary then the extension is a ${}^*$-representation.
\end{corollary}
\begin{proof}
The assumption $\Sp\rho\subseteq\Sp A$ implies $\ker\rho\supseteq
A^\bot$, where $\rho$ in the left-hand side is the extension of
$\rho$ to $M(G)$ defined by {\rm(\ref{exten})}. Hence, the
right-hand side of {\rm(\ref{exten})} is independent of the choice
of $\mu\in\cM_\vf$. Thus, the extension of $\rho$ to $\cM_A$ is
well-defined. Standard arguments show that it is strongly
continuous and, if $\rho$ is unitary, a ${}^*$-representation.
\end{proof}
\begin{proof}[Proof of Theorem~\ref{semit}] The multiplication in
$\cM_A$ is continuous due to Theorem~\ref{mabal} since the unit ball
in $\BL(A)$ with the strong operator topology is a topological semigroup.
By {\rm(\ref{convo})} and {\rm(\ref{iodef})}, $\de_g*\de_h=\de_{gh}$
for all $g,h\in G$; since $A$ is separating,
$\ev$ is one-to-one. Thus, $\ev$ is an isomorphic embedding.
For each $g\in G$, the Dirac measure $\de_g$ is the unique representing
measure for $\ev_x$ due to Lemma~\ref{peakr}.
Let $\vf,\psi\in\cM_A$, $\mu\in\cM_\vf$, $\nu\in\cM_\psi$. Since $\mu$ and
$\nu$ are positive,
$$\supp\mu*\nu=\supp\mu\cdot\supp\nu.$$
If $\vf*\psi\in G$, then $\supp\mu*\nu$ is a single point. Hence,
the same is true for $\supp\mu$ and $\supp\nu$. Thus, $\vf,\psi\in
G$. The remaining assertion was already proved in
Theorem~\ref{mabal}.
\end{proof}
For abelian $G$, there is another way to define the same semigroup
structure in $\cM_A$. In this case, $\wh G$ is the group of
characters $G\to\bbT$, which can be considered as a subgroup of
the multiplicative group $C(G)^{-1}$. Characters can be defined by
the equality
\begin{eqnarray}\label{iocha}
\iota(\chi)=\chi\otimes\chi.
\end{eqnarray}
The set $\Sp A=A\mathbin\cap\wh G$ is a semigroup in $\wh G$
(obviously, the converse is also true: if $S\subseteq\wh G$ is a
semigroup that contains identity of $\wh G$, then the closure  of
its linear span in $C(G)$ is an invariant algebra). The
restriction of $\vf\in\cM_A$ to $\Sp A$ defines a homomorphism
$\Sp A\to\ov\bbD$ that satisfies the condition $\vf({\bf 1})=1$.
The set $\Hom(\Sp A,\ov\bbD)$ of all nonzero homomorphisms $\Sp
A\to\ov\bbD$ is a topological semigroup with respect to the
pointwise multiplication and the topology of pointwise
convergence. The essential part of the following theorem was
already proved in \cite{AS}.
\begin{theorem}\label{abelc}
If $G$ is abelian, then the restriction of functionals to $\Sp
A\subset A$ defines an isomorphism $\rho:\,\cM_A\to\Hom(\Sp
A,\ov\bbD)$ of topological semigroups. For any $\vf\in\cM_A$ and
$\chi\in S$,
\begin{eqnarray}\label{invab}
\rho(\vf^*)(\chi)=\ov{\vf(\chi)},
\end{eqnarray}
where the bar denotes the complex conjugation.
\end{theorem}
\begin{proof} 
Since $S=\Sp A=\wh G\mathbin\cap A\subset A$, the topology of
pointwise convergence on $S$ is weaker than the weak topology in
$\cM_A$. Therefore, $\rho$ is continuous. Further, $\rho$ is
one-to-one because $A_\fin$ is dense in $A$. Equality
{\rm(\ref{iocha})}, together with {\rm(\ref{convo})}, implies that
$\rho$ is a homomorphism. Clearly, $\chi^\star=\chi$ for each
$\chi\in\wh G$. Hence, {\rm(\ref{invab})} follows from
{\rm(\ref{invma})}. It remains to prove that $\rho$ is surjective
(recall that $\cM_A$ is compact). Each $\vf\in\Hom(S,\ov\bbD)$ can
be extended to a continuous linear functional on $l^1(S)\subseteq
l^1(\wh G)$ since $\vf$ is bounded on $S$. Spaces $l^1(\wh G)$ and
$l^1(S)$ are convolution Banach algebras. Evidently, $\vf$ is
multiplicative on $l^1(S)$. The maximal ideal space of $l^1(\wh
G)$ is known to be equal to $G$.
By the Gelfand-Naimark formula,
$\lim_{n\to\infty}\|f^n\|^{\frac1n}_{l^1(\wh G)}=\|f\|_{C(G)}$ for
any $f\in l^1(\wh G)$.  Thus, relations
$|\vf(f)|=|\vf(f^n)|^{\frac1n}$ and
$|\vf(f^n)|\leq\|f^n\|_{l^1(\wh G)}$ imply
$|\vf(f)|\leq\|f\|_{C(G)}$. Since $l^1(S)$ can be identified with
a dense subalgebra of $A$, $\vf$ admits a continuous extension to
$A$.
\end{proof}
In the sequel, we omit $*$ in the notation for the multiplication
in $\cM_A$. The simplest (and leading) example, which illustrates
constructions of this section, is the following one.
\begin{example}\label{disca}\rm
The algebra $A(\bbD)$ of all functions that are analytic in $\bbD$
and continuous in $\ov\bbD$ can be considered as a subalgebra of
$C(\bbT)$. It is an invariant algebra on $\bbT$,
$\cM_{A(\bbD)}=\ov\bbD$, and the above multiplication in
$\cM_{A(\bbD)}$ coincides with the multiplication of complex
numbers. Indeed, for each $z\in\bbD$ the operator $R_z$ defined by
$R_zf(\ze)=f(\ze z)$, $\ze\in\bbT$, commutes with all translations
in $\bbT$ and that $R_zR_w=R_{zw}$. The involution ${}^*$ is the
complex conjugation in $\ov\bbD$ since it is the unique continuous
semigroup automorphism of $\ov\bbD$ which coincides with the
inversion on $\bbT$.
\end{example}
\begin{example}\label{uball}\rm
Analogously, the algebra $A(B_n)$ of functions analytic in the
open unit matrix ball $B_n=\{Z\in\BL(\bbC^n):\,Z^*Z<1\}$ and
continuous up to the boundary can be considered as an invariant
algebra on the group $\rU(n)$, the multiplication
{\rm(\ref{convo})} and the involution {\rm(\ref{invma})} coincides
with the standard ones. The proof is almost word for word as
above.
\end{example}
We conclude this section with a description of the centre of
$\cM_A$. Let $\Inn(G)$ denote the group of inner automorphisms
$x\to g^{-1}xg$ of $G$; they naturally extend to automorphisms of
$\cM_A$.
\begin{proposition}\label{centr}
Let $\ze\in\cM_A$. Following assertions are equivalent:
\begin{itemize}
\item[\rm1)] $\ze\in Z(\cM_A)$;
\item[\rm2)] $\ze\in Z(G)$;
\item[\rm3)] $\cM_\ze$ contains an $\Inn(G)$-invariant measure.
\end{itemize}
\end{proposition}
\begin{proof}
The implication 1)$\ \Longrightarrow\ $2) is trivial. If 2) is
true then the convex weakly closed set $\cM_\ze$ is
$\Inn(G)$-invariant. Hence it contains an $\Inn(G)$-fixed point
which is an $\Inn(G)$-invariant measure. The centre of the
convolution algebra $M(G)$ is exactly the set of such measures.
Due to Proposition~\ref{hosub}, 3) implies 1).
\end{proof}

\section{Averaging over a subgroup and restriction to a subgroup}
In this section, $H$ denotes a closed subgroup of $G$. The {\it
left averaging operator over $H$} is defined by
\begin{eqnarray}\label{aveop}
L_H=\int_HL_h\,dh, 
\end{eqnarray}
where $dh$ is the Haar measure on $H$.  It is well-defined in any
left invariant closed subspace of $C(G)$. If $h\in H$, then
$L_hL_H=L_HL_h=L_H$; moreover, $L^2_H=L_H$, and $L_H$ commutes
with $R_\vf$ for all $\vf\in \cM_A$. If $H$ is normal then $L_H$
commutes with left and right translations and coincides with the
(naturally defined) right averaging operator.
\begin{proposition}\label{averr}
The space $B=L_HA$ is a closed right invariant subalgebra of $A$. The class
$$H'=\{h\in G:\,h\eqv{B}e\}$$
includes $H$, is a closed subgroup of $G$, which is normal if $H$
is normal, and a $p$-set (a peak set, if $G$ is a Lie group).
Moreover, $L_H=L_{H'}$ in $A$ and
\begin{eqnarray}\label{averd}
B=\{f\in A:\,L_hf=f\ \text{for all}\ h\in H'\}.
\end{eqnarray}
The equality $H'=H$ holds if and only if $H$ is a $p$-set for $A$.
\end{proposition}
\begin{proof}
Since $L_H$ commutes with $R_g$ for all $g\in G$, the algebra $B$
and the equivalence ${}\eqv{B}{}$ are right invariant. Taken
together with the inclusion $e\in H'$, this implies
$H'h=H'h^{-1}=H'$ for all $h\in H'$. Thus, $H'$ is a subgroup.
Obviously, it is closed and includes $H$. If $H$ is normal, then
$L_H$ commutes with $L_g$ for all $g\in G$, hence $B$,
${}\eqv{B}{}$ are left invariant and $H'$ is normal. If $h\in H'$
and $f\in B$, then $L_hf(e)=f(e)$ by definition of $H'$; replacing
$f$ by $R_gf$, $g\in G$, we get $L_hf=f$ and, integrating over
$H'$, $L_{H'}f=f$. This proves the nontrivial part of
{\rm(\ref{averd})}. Thus, $B$ is the set of all functions in $A$
that are constant on classes $H'g$ for all $g\in G$ and can be
considered as a separating $G$-invariant function algebra on the
homogeneous space $M=H'\mathbin\backslash G$, where $G$ acts on
$M$ by right. Due to Lemma~\ref{peakr}, the class $H'$ is a
$p$-point for $B$ on $M$ (a peak point if $G$ is a Lie group);
since $B\subseteq A$, $H'$ is a $p$-set (peak set) for $A$ on $G$.
If $H$ is a $p$-set and $g\notin H$ then there exists a peak set
$P$ such that $H\subseteq P$ but $g\notin P$.  Let $f$ be a peak
function for $P$. Then $L_Hf$ is a peak function for some peak set
$P'$ such that $H\subseteq P'\subseteq P$. Hence,
$L_Hf(e)=L_Hf(h)$ implies $h\in P'$; therefore, $g\notin H'$.
Thus, $H=H'$ if $H$ is a $p$-set. Remaining assertions are clear.
\end{proof}
The equivalence $\eqv{B}$ naturally extends to $\cM_A$:
$\vf\eqv{B}\psi$, if $f(\vf)=f(\psi)$ for all $f\in B$.
Factorizing by $\eqv{B}$, we get the image of $\cM_A$ under the
mapping $\rho:\,A'\to B'$ dual to the embedding $B\to A$. Clearly,
$\rho(\cM_A)\subseteq\cM_B$. In general, $\rho(\cM_A)\neq\cM_B$
but in the setting above $\rho$ is surjective due to the following
(obvious) equality:
\begin{eqnarray}\label{averi}
L_H(ab)=(L_Ha)b \quad\text{for all}\ a\in A,\ b\in B.
\end{eqnarray}
\begin{proposition}\label{surme}
The above mapping $\rho:\,\cM_A\to\cM_B$ is surjective.
\end{proposition}
\begin{proof}
Let $I$ be a proper ideal in $B$. If $I$ does not generate a proper ideal
in $A$ then there exist functions $b_1,\dots b_n\in I$ and
$a_1,\dots,a_n\in A$ such that $a_1b_1+\dots+a_nb_n=1$. Then,
by {\rm(\ref{averi})},
\begin{eqnarray*}
1=(L_Ha_1)b_1+\dots+(L_Ha_n)b_n\in I
\end{eqnarray*}
contradictory to the assumption that
$I$ is proper.  Hence each maximal proper ideal of $B$ is contained
in some maximal proper ideal of $A$.
\end{proof}
Let $M=H\mathbin\backslash G$ be a right homogeneous space of $G$.
In the following theorem, the algebra $A^H=L_HA$ is considered as
a subalgebra of $C(M)$ as well as a subalgebra of $C(G)$.
\begin{theorem}\label{surje}
The action of $G$ on $M$ extends to the action of $\cM_A$ on $\cM_{A^H}$ by
\begin{eqnarray}\label{actmb}
f(\psi\vf)=R_\vf f(\psi)
\end{eqnarray}
for all $f\in A^H$, $\psi\in\cM_{A^H}$, where $\vf\in\cM_A$ and the
right side defines the left one.
For each $p\in M$,
\begin{eqnarray}\label{suror}
p\cM_A=\cM_{A^H}.
\end{eqnarray}
\end{theorem}
\begin{proof}
Obviously, $A^H$ is right invariant.
Since ${A^H}\subseteq A$, each irreducible component of the right regular
representation $R$ of $G$ in $A^H$ is contained in $\Sp A$.
By {\rm(\ref{intrf})}, $R$ extends to the representation of $\cM_A$ in $A^H$;
clearly, this is equivalent to {\rm(\ref{actmb})} and defines
the action of $\cM_A$ on $\cM_{A^H}$ that commutes with the mapping
$\rho$ above.
Since $g\cM_A=\cM_A$ for any $g\in G$, {\rm(\ref{suror})}
follows from Proposition~\ref{surme}.
\end{proof}
\begin{proposition}\label{resta}
If $H$ is a closed subgroup of $G$, then $\wh H$ is a closed
subsemigroup of $\cM_A$. Moreover, the closure of the restriction
$A\big|_H$ is an invariant algebra $C$ on $H$, and the mapping
dual to the natural homomorphism $A\to C$ is a topological
isomorphism of ${}^*$-semigroups $\cM_C$ and $\wh H$.
\end{proposition}
\begin{proof}
According to Proposition~\ref{hosub}, the multiplication in
$\cM_A$ and in $\cM_C$ is induced by the convolution of
representing measures. If $\vf,\psi\in\wh H$, then, by
Lemma~\ref{rpeak}, there exist measures $\mu\in\cM_\vf$ and
$\nu\in\cM_\psi$ concentrated in $H$. Then
$$\supp\mu*\nu\subseteq H;$$
hence $\vf\psi\in\wh H$ by the same lemma.
Since $\cM_C$ is compact,
remaining assertions follow from this fact and the natural identification
$\cM_C=\wh H$
(which holds for any closed subset $H\subseteq G$).
\end{proof}
\begin{corollary}
The polynomially convex hull of a compact linear group $G\subseteq\GL(n,\bbC)$
is a semigroup.
\end{corollary}
\begin{proof}
We may assume $G\subseteq\rU(n)$; then Proposition~\ref{resta}
can be applied to the algebra of Example~\ref{uball}.
\end{proof}

By the well-known structure theorem for compact groups, each
neighbourhood of $e$ contains a closed normal subgroup $H$ such
that $G/H$ is a Lie group. For any pair $H_\al\subseteq H_\be$ of
normal subgroups there is the natural homomorphism
$\vf_{\al\be}:\,G_\al\to G_\be$, where $G_\al=G/H_\al$ and
$G_\be=G/H_\be$. Let $\frH=\{H_\al\}_{\al\in I}$, where $I$ is a
set of indices, be a family of closed normal subgroups in $G$ such
that
\begin{eqnarray} 
&H,P\in\frH\quad\Rightarrow H\mathbin\cap P\in\frH,
\label{iner1}\\ 
&\bigcap_{H\in\frH}H=\{e\}\label{iner2}.
\end{eqnarray}
Then $G$ is the inverse limit of Lie groups $G/H$ over $\frH$,
i.e. it is topologically isomorphic
to the subgroup of the topological Cartesian product of these groups
consisting of families $(g_\al)_{\al\in I}$
such that $g_\be=\vf_{\al\be}(g_\al)$ if $H_\al\subseteq H_\be$.
By the following theorem, the semigroup $\cM_A$ can be realized analogously.
Let $\frL=\frL(G)$ be the family of all normal subgroups $H\subseteq G$
such that $G/H$ is a Lie group and
let $\frP=\frP(G,A)$ be its subfamily consisting of those $H$ which are
peak sets for $A$.
\begin{theorem}\label{invli}
Let $A$ be a separating invariant algebra on a compact Hausdorff
group $G$. For any $H\in\frP$, $A^H=L_HA$ is a separating
invariant algebra on the Lie group $G/H$.  The projection
$\rho_H:\,\cM_A\to\cM_{A^H}$ dual to the embedding $A^H\to A$ is a
surjective continuous homomorphism of semigroups. The family
$\frP$ satisfies {\rm(\ref{iner1})} and {\rm(\ref{iner2})}; as a
compact topological semigroup, $\cM_A$ is the inverse limit of
$\cM_{A^H}$, where $H$ runs over $\frP$.
\end{theorem}
\begin{proof}
Since $H$ is normal, its Haar measure $m_H$ belongs to the centre
of the convolution algebra $M(G)$. The projection
$\cM_A\to\cM_{A^H}$ is a homomorphism because the multiplication
in $\cM_A$ is induced by the convolution and the measure $m_H$ is
convolution idempotent. Obviously, it is continuous. By
Theorem~\ref{surje}, it is surjective. Let $u$ be the peak
function for $H\in\frH$. Then $L_Hu$ is the peak function for the
identity element of $G/H$. Hence, the equivalence $\eqv{A^H}$ on
$G/H$ is trivial, i.e., $A^H$ is separating on $G/H$. If
$H,P\in\frP$ and $u,v$ are peak functions for $H$ and $P$,
respectively, then  $uv$ has a peak on $H\mathbin\cap P$. This
proves {\rm(\ref{iner1})}. The family of operators $L_H$,
$H\in\frL$, is an approximate identity in $C(G)$, hence in $A$.
Therefore, if $g\in G$ and $f(g)\neq f(e)$ for some $f\in A$, then
$L_Hf(g)\neq L_Hf(e)$ for some $H$ as above. Taken together with
Proposition~\ref{averr}, particularly {\rm(\ref{averd})}, this
proves {\rm(\ref{iner2})}. The remaining assertion is clear.
\end{proof}
The theorem above reduces most of the problems on invariant
algebras to the case of Lie groups. To illustrate this, we prove a
lemma which will be used in characterization of compact groups $G$
which have the following property:
\begin{itemize}
\item[(\sf{SW})] each invariant algebra on $G$ is self-adjoint.
\end{itemize}
We need the following result, which is a consequence of
\cite[Theorem 6.1]{Bj}.
\begin{theorem}\label{a=cg;}
Let an invariant algebra $A$ be separating on a compact group $G$
and let $\cM_A=G$. Then $A=C(G)$.
\end{theorem}
For compact Lie groups, this theorem easily follows from known
facts on polynomial approximation on polynomially convex
submanifolds in $\bbC^n$.
\begin{lemma}\label{redsw}
The condition {\rm(\sf{SW})} holds for $G$ if and only if
{\rm(\sf{SW})} is true for all groups $G/H$, where $H$ runs over
$\frL$.
\end{lemma}
\begin{proof}
Let $A$ be an invariant algebra on $G$ which is not self-adjoint.
We may assume that $A$ is separating. By Theorem~\ref{a=cg;},
$\cM_A\neq G$; due to Theorem~\ref{invli}, $\cM_{A^H}\neq G/H$ for some
$H\in\frP$. This proves the nontrivial part of the lemma.
\end{proof}

\section{Antisymmetric algebras and idempotents}\label{antii}

Let $\cI_A$ be the set of all idempotents in $\cM_A$. For
$j\in\cI_A$, put
\begin{eqnarray*}
&G_j=\{g\in G:\,gj=jg=j\},\quad N_j=Z_G(j)=\{g\in G:\,gj=jg\},\\
&G_j^l=\{g\in G:\,gj=j\},\quad
G_j^r=\{g\in G:\,jg=j\},\\
&G^j=N_jj\cong N_j/G_j,\quad \cM^j=\{\vf\in\cM_A:\,\vf j=j\vf=\vf\}=j\cM_Aj,
\end{eqnarray*}
and denote by $m_j$ the normalized Haar measure of $G_j$. There is
a natural order in the set of idempotents of any semigroup: $j\leq
k$ if $jk=kj=j$.

The following proposition will be often used in this section.
Hence references of the type ``Proposition~\ref{idems}, x)" will
be abbreviated to ``x)".
\begin{proposition}\label{idems}
Let $j\in\cI_A$. Then
\begin{itemize}
\item[\rm1)] $j^*=j$;
\item[\rm2)] $G_j=G_j^l=G_j^r$;
\item[\rm3)] $G_j$ is a $p$-set for $A$;
\item[\rm4)] $m_j\in{\cal M}_j$;
\item[\rm5)] for any $k\in\cI_A$, $j\leq k$ is equivalent to
$G_j\supseteq G_k$;
\item[\rm6)] $N_j=N_G(G_j)$; in particular, $j\in Z(\cM_A)$
if and only if $G_j$ is normal;
\item[\rm7)] $\vf,\psi\in\cM^j$ and $\vf\psi=j$ imply $\vf,\psi\in G^j$.
\end{itemize}
\end{proposition}
\begin{proof}
Since $R_j$ is an endomorphism of $A$ commuting with left
translations and $R_j^2=R_j$, it is a projection onto a closed
left invariant subalgebra $B\subseteq A$. The left orbit ${\cal
O}_j =\{gj:\,g\in G\}$ in $\cM_A$ is the homogeneous space
$G/{G_j^l}$; $A|_{{\cal O}_j}$ is a separating subalgebra $\tilde
B$ of $C({\cal O}_j)$. Due to the identity $f(gj)=R_j f(g)$, which
holds for all $f\in A$ and $g\in G$, $\tilde B$ may be identified
with $B$. According to Lemma~\ref{peakr}, each point of ${\cal
O}_j$ is a $p$-point; hence, ${G_j^l}$ is a $p$-set for $A$. By
Lemma~\ref{rpeak},
\begin{eqnarray}\label{suppm}
\supp\mu\subseteq{G_j^l}
\end{eqnarray}
for any $\mu\in\cM_j$. Since the set ${\cal M}_j$ is  weakly
compact, convex and $G_j^l$-invariant by left this implies that
${\cal M}_j$ contains the Haar measure $m_j^l$ of $G_j^l$. The
measure $m_j^l$ is inversion invariant and positive. Hence
$(m_j^l)^*=m_j^l$;  this proves 1). The equality $gj=(j g^{-1})^*$
and 1) imply $G_j^l=G_j^r$, consequently  2). Since ${G_j^l}$ is a
$p$-set for $A$ and $m_j^l\in\cM_j$, 2) implies  3) and 4). The
inclusion $G_j\supseteq G_k$ is equivalent to
$m_j*m_k=m_k*m_j=m_j$; this proves 5). If $g\in N_G(G_j)$, then
$\de_g*m_j=m_j*\de_g$, whence $gj=jg$; conversely, if
$gjg^{-1}=j$, then
$$\supp(\de_g*m_j*\de_{g^{-1}})=gG_jg^{-1}\subseteq G_j$$
by {\rm(\ref{suppm})}. Thus, $g\in N_G(G_j)$. It remains to prove
7). If $\mu\in\cM_\vf$, $\nu\in\cM_\psi$, then
$m_j*\mu*m_j\in\cM_\vf$ and $m_j*\nu*m_j\in\cM_\psi$. Hence, we
may assume that $m_j*\mu*m_j=\mu$, $m_j*\nu*m_j=\nu$. Then
$\mu*\nu=m_j$ by {\rm(\ref{suppm})}. Set $P=\supp\mu$,
$Q=\supp\nu$. By equalities above,
\begin{eqnarray*}
G_jPG_j=PG_j=G_jP=P,\quad G_jQG_j=QG_j=G_jQ=Q,\quad PQ=G_j.
\end{eqnarray*}
It follows that $P=gG_j=G_jg$, $Q=g^{-1}G_j=G_jg^{-1}$ for any
$g\in P$. Clearly, $g\in N_G(G_j)$. Further,
$m_j*\mu*\de_{g^{-1}}=m_j$, whence $j\vf g^{-1}=j$. Since $\vf\in
j\cM_Aj$, this implies $\vf=gj=jg\in G^j$. Analogously,
$\psi=g^{-1}j=jg^{-1}\in G^j$.
\end{proof}
\begin{lemma}\label{lowid}
A closed subsemigroup $\cS\subseteq\cM_A$ generated by a subset
$\cJ\subseteq\cI$ is a semigroup with zero.
\end{lemma}
\begin{proof}
By Theorem~\ref{mabal}, the mapping $\vf\to R_\vf$ identifies
$\cM_A$ with the semigroup $\frM_A\subseteq\BL(H^2)$, which acts
in $H^2$. Hence, $\frS=\{R_s:\,s\in\cS\}$ is a strongly closed
semigroup acting in $H^2$, and $\frS^*=\frS$. By 1) and
Theorem~\ref{mabal}, if $j\in\cI$, then $R_j$ is an orthogonal
projection. Set
\begin{eqnarray*}
L=\bigcap_{j\in \cJ}R_jH^2
\end{eqnarray*}
and let $E$ be the orthogonal projection onto $L$. Clearly, $L$
and $L^\bot$ are left invariant and $\frS$-invariant; moreover,
$ER_s=R_sE=E$ for all $s\in\frS$. It is sufficient to prove that
$E\in\frS$. Otherwise, there exists a neighbourhood (in the strong
operator topology) of $E$ in $\BL(H^2)$ which separates $\frS$ and
$E$. Then there exist $f_1,\dots,f_k\in L^\bot$ such that
\begin{eqnarray}\label{sepne}
\|Sf_l\|_{H^2}>1\quad\mbox{for all}\quad
S\in\frS\quad\mbox{and}\quad l=1,\dots,k.
\end{eqnarray}
For $\tau\in\wh G$, set $L^\bot_\tau=L^\bot\mathbin\cap M_\tau$.
These spaces are left invariant and $\frS$-invariant, and pairwise
orthogonal. Moreover, their algebraic sum $\sum_{\tau\in\wh
G}L^\bot_\tau$ is  dense in $L^\bot$. Since $\frS$ is strongly
compact, we may assume that the functions $f_1,\dots,f_k$, are
contained in a finite sum $N$ of these spaces. Since $L\cap
N=\{0\}$ and $N$ is finite dimensional, there exists a finite set
$J=\{j_1,\dots,j_m\}\subseteq\cI$ such that $N\bigcap_{j\in
J}R_jH^2=\{0\}$. Then, for each $f\in N$, there exists $j\in J$
such that $\|R_jf\|<\|f\|$. It follows that the eigenvalues of the
self-adjoint operator
\begin{eqnarray*}
&R_s=R_{j_1}R_{j_2}\dots R_{j_m}R_{j_m}\dots R_{j_2} R_{j_1},
\quad s=(j_m\dots j_2j_1)^*(j_m\dots j_2j_1),
\end{eqnarray*}
in $N$ are strictly less than $1$. Therefore,
$\|R_{s^n}f_l\|_{H^2}\to0$ as $n\to\infty$ for each $l$,
contradictory to (\ref{sepne}).
\end{proof}
An ordered set is called a {\it complete lattice} if each its bounded
subset $I$ has the least upper bound $\sup I$ and the greatest lower bound
$\inf I$.
\begin{theorem}\label{icola}
The set $\cI_A$ is a complete lattice; moreover, $\sup I$ and
$\inf I$ are well-defined for each subset $I\subseteq\cI_A$. In
particular, $\cI_A$ contains the least idempotent $\eps$ which
belongs to $Z(\cM_A)$.
\end{theorem}
\begin{proof}
If $jk=kj=j$ for all $k\in I$ then $j\vf=\vf j=j$ for all $\vf$ in
the semigroup generated by $I$. Clearly, the zero of this
semigroup, which exists by  Lemma~\ref{lowid}, is $\inf I$.
Further, $\sup I=\inf J$, where $J$ is the set of all upper bounds
for $I$ (note that $J\neq\emptyset$ because $e$ is the greatest
element of $\cI_A$ and that $\vf k=k\vf=k$ for all $k\in I$ and
$\vf$ in the closed semigroup generated by $J$). By 5),
$G_\eps\supseteq G_j$ for all $j\in\cI_A$; in particular
$g^{-1}G_\eps g\subseteq G_\eps$ for all $g\in G$. Hence, $G_\eps$
is normal. By 4) and Proposition~\ref{centr}, $\eps\in Z(\cM_A)$.
\end{proof}

A function algebra $B\subseteq C(Q)$ is called {\it antisymmetric}
if each real valued function $f\in B$ is constant; a set $E$ is
called a {\it set of antisymmetry} if any function $f\in B$, which
is  real valued on $E$, is constant on $E$. Clearly, if two sets
of antisymmetry intersects, then their union is a set of
antisymmetry, and any set of antisymmetry is contained in a
maximal one. By the Shilov--Bishop decomposition theorem, a
continuous function $f$ is contained in $B$ if and only if
$f|_E\in B|_E$ for every maximal set of antisymmetry $E$. Let
$B_\bbR$ be the algebra of all real valued functions in $B$.
According to the Stone--Weierstrass theorem, $B_\bbR$ can be
identified with the algebra of all real valued continuous
functions on $Q/\eqv{B_\bbR}$. In general, classes of
${\,}\eqv{B_\bbR}{\,}$ are not sets of antisymmetry, but for
invariant algebras they are.
\begin{theorem}\label{antis}
The set $G_\eps$, where $\eps$ is the least idempotent, is the
maximal set of antisymmetry for $A$. It is a normal subgroup, all
other maximal sets of antisymmetry are cosets $gG_\eps$, $g\in G$,
and $A_{\bbR}$ consists of all real valued continuous functions
which are constant on them. Furthermore, following assertions are
equivalent:
\begin{itemize}
\item[\rm a)] $A$ is antisymmetric;
\item[\rm b)] $\cM_A$ is a semigroup with zero;
\item[\rm c)] the Haar measure of $G$ is multiplicative on $A$.
\end{itemize}
\end{theorem}
\begin{proof}
The subgroup $G_\eps$ is normal by Theorem~\ref{icola} and 6). The
set $\supp\mu$ for any $\mu\in\cM_\vf$ and each $\vf\in\cM_A$ is
known to be the set of antisymmetry (indeed, if $f\in A$ is real
nonconstant on $\supp\mu$, then $\vf\left((f-f(\vf))^2\right)>0$,
but this cannot be true since $\vf\left(f-\vf(f)\right)=0$ and
$\vf$ is multiplicative). According to 4), $G_\eps$ is a set of
antisymmetry. By Theorem~\ref{invli} and Theorem~\ref{a=cg;},
$A^{G_\eps}=C(G/G_\eps)$. Hence, $G_\eps$ is a maximal set of
antisymmetry and the same is true for its shifts $gG_\eps$. It
follows that $A$ is antisymmetric if and only if $G_\eps=G$, i.e.
if b) is true; b) implies c) by 4). Since the support of a
multiplicative measure is a set of antisymmetry, c) implies a).
\end{proof}
For abelian $G$, there is another simple criterion of
antisymmetry, in terms of the semigroup $S=\Sp A$. Precisely, $A$
is antisymmetric if and only if
\begin{eqnarray}\label{antab}
S\mathbin\cap(-S)=\{0\}
\end{eqnarray}
(in the additive notation). Indeed, if $\chi,\ov\chi\in A$, then
$\Re\chi\in A$. Conversely, $A$ contains all characters in the
support of the Fourier transform of any $f\in A$, but for real $f$
the support is invariant with respect to the inversion in $\wh G$;
hence, (\ref{antab}) cannot be true if $f$ is nonconstant.

\section{Polar decomposition, one parameter
semigroups, and analytic structure} This section overlaps with the
paper \cite{Gi79}.
\begin{lemma}\label{l2uni}
Let $B\subseteq C(Q)$ be a function algebra (in general,
nonclosed) on compact $Q$ and $\mu\in M(Q)$ be a probability
measure on $Q$. Suppose $\supp\mu=Q$ and let $E$ be an
endomorphism of $B$ which is bounded as a linear operator $B\to B$
with respect to the norm in $L^2(Q,\mu)$. Then $E$ is bounded in
$B$ with respect to the $\sup$-norm; moreover,
$\|E\|_{\BL(B)}\leq1$.
\end{lemma}
\begin{proof}
Suppose that $f\in B$ satisfies inequalities $\|f\|<c<1$ and
$\|Ef\|>1$. Then the set $U=\{q\in Q:\, |f(q)|>1\}$ is open and
nonvoid. Hence $\mu(U)>0$. Therefore,
\begin{eqnarray*}
&\|Ef^n\|_{L^2(Q,\mu)}=\|(Ef)^n\|_{L^2(Q,\mu)}\geq\sqrt{\mu(U)},\\
&\|f^n\|_{L^2(Q,\mu)}\leq\|f^n\|=\|f\|^n\leq
c^n\to0\quad\text{as}\ n\to\infty,
\end{eqnarray*}
contradictory to the assumption.
Thus $\|E\|_{\BL(B)}\leq1$.
\end{proof}
Let $B\subseteq C(Q)$ be a subalgebra that admits an orthogonal
grading by a commutative semigroup $\La$ (in the multiplicative
notation), i.e.,
\begin{eqnarray}
&B=\sum_{r\in \La}B_r,\label{sumbs}\\
&B_t\cdot B_r\subseteq B_{tr}\quad\text{for all}\quad t,r\in \La,\label{probs}\\
&B_t\perp B_r\quad \text{in} \quad  L^2(Q,\mu)\quad \text{if}\quad
t\neq r,\label{perbs}
\end{eqnarray}
where the sum in (\ref{sumbs}) is algebraic. 
We do not assume $B_t\neq0$.
\begin{lemma}\label{chibo}
Let $B$ be as above, $\chi:\,\La\to\ov\bbD$ be a homomorphism of
semigroups, and let a linear operator $E:\,B\to B$ be defined by
$Ef=\chi(t)f$ for $f\in B_t$, $t\in\La$. Then $E$ is bounded in
$B$ with respect to the $\sup$-norm; moreover,
$\|E\|_{\BL(B)}\leq1$.
\end{lemma}
\begin{proof}
A calculation shows that $E$ is an endomorphism of $B$.
Since $\chi$ is bounded and $B_t$ are pairwise orthogonal,
$E$ is bounded in $B$ with respect to the norm of $L^2(Q,\mu)$,
and the assertion follows from Lemma~\ref{l2uni}.
\end{proof}
Set
\begin{eqnarray}\label{defsa}
\cS_A=\{s\in\cM_A:\,s^*=s,\ R_s\geq0\ \text{in}\ H^2\}.
\end{eqnarray}
We denote by $\bbC^+$ the closed right halfplane in $\bbC$,
$\bbR^+=[0,\infty)$. For any $\vf\in\cM_A$, $\vf^*\vf\in \cS_A$
due to Theorem~\ref{mabal}. Let $\ga:\,\bbR^+\to\cM_A$  be a
continuous homomorphism (i.e., a one parameter semigroup). We say
that $\ga$ is a {\it ray} if it is symmetric: $\ga(t)^*=\ga(t)$
for all $t\geq0$. Then $\ga(\bbR^+)\subseteq \cS_A$. Let $\frR$ be
the family of all rays in $\cM_A$. Note that $\ga(0)\in\cI_A$. For
$k\in\cI_A$, set
\begin{eqnarray}
&\frR^k=\{\ga\in\frR:\,\ga(0)=k\},\label{frrk.}  \\
&\cS_A^k=\mathop\bigcup\nolimits_{\ga\in\frR^k}{\ga(\bbR^+)}\label{sakk.}.
\end{eqnarray}
A continuous homomorphism $\ga:\,\bbC^+\to\cM_A$ will be called a
{\it complex ray} if $f(\ga(z))$ is analytic in the open right half-plane
for all $f\in A$.
\begin{theorem}\label{rays.}
For any $s\in\cS_A$, there exists the unique ray $\ga$ such that
$\ga(1)=s$. Each ray admits the unique extension to a complex ray.
For any complex ray $\ga$, there exists the limit
\begin{eqnarray}\label{infty}
\lim_{\Re z\to\infty}\ga(z)=:\ga(\infty)\in\cI_A.
\end{eqnarray}
Furthermore, $\ga(it)$, $t\in\bbR$, is one parameter group
in $G^j$, where $j=\ga(0)$.
\end{theorem}
\begin{proof}
Let $\La$ be the semigroup generated by the eigenvalues of $R_s$.
Clearly, $\La\subseteq[0,1]$. For $t\in\La$, let $A_t$ be the
$t$-eigenspace of $R_s$ in $A_\fin$ (if $t$ is not an eigenvalue,
then $A_t=0$). We claim that these spaces satisfies assumptions of
Lemma~\ref{chibo} with $A_\fin$ as $B$. Since $R_s$ is
self-adjoint, $A_t\perp A_r$ in $L^2(G)$ if $t,r\in\La$ and $t\neq
r$; hence {\rm(\ref{perbs})} holds. Any minimal bi-invariant
subspace of $H^2$ is $R_s$-invariant and finite dimensional; hence
$A_\fin=\sum_{t\in\La}A_t$. This proves {\rm(\ref{sumbs})}. If
$p\in A_t$, $q\in A_r$, then $R_s(pq)=(R_sp)(R_sq)=trpq$, i.e.
$pq\in A_{tr}$. Thus {\rm(\ref{probs})} is true. Put
\begin{eqnarray}\label{chien}
\chi_z(t)=\left\{\begin{array}{cc}e^{z\log t},&1\geq t>0,\cr
                               0,&t=0,\end{array}\right.
\end{eqnarray}
where $\log t$ is real. Then $\chi_z$ is a bounded homomorphism
$\La\to\bbD$ for any $z\in\bbC^+$, and we may apply
Lemma~\ref{chibo} to $A_\fin$.  According to it, the
multiplication by $\chi_z(t)$ on $A_t$ extends to a continuous
endomorphism $E:\,A\to A$. It commutes with left translations
because $R_s$ has this property. Hence $E=R_{\ga(z)}$ for some
$\ga(z)\in\cM_A$. Since $\chi_{z+w}=\chi_z\chi_w$, $z\to
R_{\ga(z)}$ is a homomorphism $\bbC^+\to\frM_A$ and
$\ga:\,\bbC^+\to\cM_A$ is a homomorphism by Theorem~\ref{mabal}.
If $f\in A_\fin$, then $f(\ga(z))$ is analytic inside $\bbC^+$
being the finite sum of the type $\sum c_ke^{\mu_kz}$. The mapping
$z\to\chi_z(t)$ is continuous in $\bbC^+$ for all $t\in\La$ and
bounded, hence, the homomorphism $z\to R_{\ga(z)}$ is continuous
with respect to the strong operator topology in $\BL(H^2)$. Due to
Theorem~\ref{mabal}, $z\to\ga(z)$ is continuous on $\bbC^+$. Since
$A_\fin$ is dense in $A$, $f(\ga(z))$ is analytic in the open
half-plane and continuous in $\bbC^+$ for all $f\in A$. The
uniqueness of the ray follows from the uniqueness of fractional
nonnegative powers of a nonnegative operator in a Hilbert space.
The analytic extension is evidently unique. Equality
{\rm(\ref{infty})} is a consequence of {\rm(\ref{chien})}:
$R_{\ga(\infty)}$ is the projection to the closure of $A_1$. The
last assertion follows from Proposition~\ref{idems}, 7) due to
equalities $j\ga(it)=\ga(it)j=\ga(it)$ and $\ga(it)\ga(-it)=j$,
where $t\in\bbR$.
\end{proof}
We shall also denote $\ga(z)$ by $s^z$. Clearly, $s^0,s^\infty\in\cI_A$.
\begin{corollary}\label{comst}
For any $s\in\cS_A$ and  $z\in\bbC^+$, $Z_G(s)\subseteq Z_G(s^z)$.
Moreover, if $t>0$, then $Z_G(s)=Z_G(s^t)$.
\end{corollary}
\begin{proof}
Suppose that $s$ is a fixed point of the automorphism $\vf\to
g^{-1}\vf g$. Then $s^z$, $\Re z>0$, also is a fixed point due to
the uniqueness part of Theorem~\ref{rays.}, and for $\Re z\geq0$
by the continuity. Therefore, $Z_G(s)\subseteq Z_G(s^z)$. If
$t>0$, then $s^t\in\cS_A$ and the inverse inclusion holds since
$s=(s^t)^{\frac1t}$.
\end{proof}
In the proof of the following theorem, we refer to
``Proposition~\ref{idems},~x)'' as ``x)''.
\begin{theorem}[\rm Polar decomposition]\label{polar}
For any $\vf\in\cM_A$, there exist $s\in\cS_A$ and $g\in G$ such
that $\vf=gs$. In this decomposition, $s$ is determined by $\vf$
uniquely and $gs=hs$ for $h\in G$ if and only if $gG_k=hG_k$,
where $k={s^0}$; in particular, if $k=e$, then the decomposition
is unique.
\end{theorem}
\begin{proof}
Let $R_\vf=US$ be the polar decomposition of $R_\vf$ in the
Hilbert space $H^2$. Then
$S^2=R_{\vf}^*R_\vf=R_{\vf^*}R_\vf=R_{\vf^*\vf}$ by
Theorem~\ref{mabal}; hence, $S=R_s$ for $s=(\vf^*\vf)^{\frac12}$.
Since $R_g$ is unitary for any $g\in G$, this proves the
uniqueness of the component $s\in\cS_A$ in the polar decomposition
$\vf=gs$, if it exists.

Suppose $s=k$. Then $k=k^2=\vf^*\vf$. Pick $\mu\in\cM_\vf$,
$g\in\supp\mu$, and set $\ka=g^{-1}\vf$ (hence $\vf=g\ka$),
$\nu=\de_{g^{-1}}*\mu\in\cM_\ka$.
Then $\ka^*\ka=k$ and $e\in\supp\nu$;  this yields, respectively,
\begin{eqnarray*}
&\supp(\nu^**\nu)\subseteq G_k,\\
&\supp\nu\subseteq\supp\nu^*\cdot\supp\nu=\supp(\nu^**\nu),
\end{eqnarray*}
where the first inclusion holds by 3) and  Lemma~\ref{rpeak} since
$\nu^**\nu$ is a representing measure for $\ka^*\ka$. Therefore,
$\supp\nu\subseteq G_k$. Due to 4), this implies $\ka k=k$.
Further, $R_\vf=UR_k$ for some unitary operator $U$; since
$k^2=k$, we have $\vf k=\vf$ and $\vf=\vf k=g\ka k=gk$. Thus, the
decomposition exists. By 4), we may assume $\supp\mu=g\supp
m_k=gG_k$. Hence, the arguments above can be applied to each $h\in
gG_k$; consequently, $\vf=gk=hk$ if  $hG_k=gG_k$. If $gk=hk$, then
$h^{-1}g\in G_k$ by 2).

The consideration above shows that the theorem holds for
$s=(\vf^*\vf)^{\frac12}\in\cI_A$. In general, the operator
$R_s^{-1}$ is well-defined on the space $A_\fin\mathbin\cap R_kA$
and the operator $E=R_\vf R^{-1}_sR_k=UR_k$ admits continuous
extension from $A_\fin$ to $H^2$. Clearly, $E$ is an endomorphism
of $A_\fin$. By Lemma~\ref{l2uni}, $E$ is bounded with respect to
the $\sup$-norm. Since $E$ commutes with left translations,
$E=R_\psi$ for some $\psi\in\cM_A$ according to
Theorem~\ref{mabal}. It follows that $\psi^*\psi=k$ and $\psi
s=\vf$. Due to the first equality, we may apply the proven
assertion. Hence, there exists $g\in G$ such that $\psi=gk$. Then
$\vf=\psi s=gks=gs$.

Evidently, $gs=hs$ if and only if $R_g=R_h$ on $R_s H^2$; since
$R_s H^2$ is a dense subspace of $R_kH^2$, this is equivalent to $gk=hk$.
The reference to 2) concludes the proof.
\end{proof}
If $G$ is abelian, then $\cM_A=\Hom(\Sp A,\bbD)$ by
Theorem~\ref{abelc}, and the polar decomposition of
$\chi\in\Hom(\Sp A,\bbD)$ is the natural polar decomposition of
functions: $\chi(x)=\rho(x)|\chi(x)|$, where $\rho\in\Hom(\Sp
A,\bbT)$. This equality uniquely defines $\rho(x)$ if
$\chi(x)\neq0$. In general, $\rho$ must be extended to $\Sp A$;
the existence of the extension follows from Theorem~\ref{polar}.

\section{Invariant algebras on Lie groups}
Everywhere in this section, $G$ is a Lie group. We keep the
notation of Section~\ref{antii} and add new: for each $j\in\cI_A$,
$Z_j$ is the identity component of the centre of $G_j$;
$\frg_j,\frn_j,\frg^j,\frz_j$ are Lie algebras of groups
$G_j,N_j,G^j,Z_j$, respectively;  $\pi_j$ is the canonical
homomorphism $N_j\to G^j=N_j/G_j$. It induces projections in
various linear spaces, which will also be denoted by $\pi_j$.
Further, $T$ denotes a maximal torus in $G$. If $G$ is connected,
then $T$ is a maximal abelian subgroup of $G$.

Each $\tau\in\Sp A$ can be extended to a representation of $\cM_A$
by Corollary~\ref{reext}. For any $j\in\cI_A$, this defines an
extension of $\tau$ to the Lie algebra $\frg^j$ by
\begin{eqnarray*}
\tau(\xi)=\frac{d}{dt}\tau(\exp(t\xi)j)\Big|_{t=0}=\tau(j)\tau(\td\xi)
=\tau(\td\xi)\tau(j),
\end{eqnarray*}
where $\xi=\pi_j\td\xi\in\frg^j$, $\td\xi\in\frn_j$. The
representation of $G^j$ and $\frg^j$ is realized in the space
$V^j_\tau=\tau(j)V_\tau$.
\begin{theorem}\label{Liera}
Let $\ga\in\frR$ and $j=\ga(0)$. There exists the unique $\xi_\ga\in\frg^j$
such that
\begin{eqnarray}\label{onepa}
\ga(it)=\exp(t\xi_\ga)=j\exp(t\xi)=\exp(t\xi)j,
\quad t\in\bbR,
\end{eqnarray}
for any $\xi\in\pi_j^{-1}\xi_\ga\subseteq\frn_j$.
The set
\begin{eqnarray}\label{coneg}
C^j=\{\xi_\ga:\,\ga\in\frR^j\}
\end{eqnarray}
is a closed convex pointed $\Ad(G^j)$-invariant cone
$C^j\subset\frg^j$, which can also be defined by
\begin{eqnarray}\label{conex}
C^j=\{\xi\in\frg^j:\,i\tau(\xi)\leq0\ \text{for all}\ \tau\in\Sp A\}.
\end{eqnarray}
Mappings $\xi_\ga\mapsto\ga$ and
$\xi_\ga\mapsto\ga(1)$ are one-to-one
on $C^j$; they identify $C^j$ with $\frR^j$ and $S_A^j$,
respectively, where $\frR^j,S_A^j$ are
defined by {\rm(\ref{frrk.})}, {\rm(\ref{sakk.})}.
\end{theorem}
\begin{proof}
The first assertion is a consequence of Theorem~\ref{rays.} (note
that the mapping $g\to jg$ on $N_j$ is the homomorphism $\pi_j$).
Clearly, {\rm(\ref{conex})} defines a convex, closed and
$\Ad(G^j)$-invariant cone. If $\xi\in C^j$, $\tau\in\Sp A$, then
$\exp(iz\tau(\xi))$ is bounded in $V_\tau$ uniformly on
$z\in\bbC^+$. Therefore, the operator group $R_{\exp(t\xi)j}$,
$t\in\bbR$, admits analytic extension to $i\bbC^+$ in $R_jM_\tau$
for all $\tau\in\Sp A$, hence to $R_jA_\fin$, and this extension
is bounded with respect to $L^2$-norm. Furthermore, the extension
is an endomorphism of the algebra $A_j=R_jA_\fin$ since
$R_{\exp(z\xi)j}\in\Hom(A_j,A_j)$ for $z\in\bbR$ and the analytic
extension agrees with the multiplication of functions. This
defines $R_{\exp(it\xi)j}$ for $t\geq0$. It is a semigroup of
endomorphisms of $A_\fin$, which commutes with left translations
and is symmetric in all $M_\tau$, $\tau\in\Sp A$. Hence,  it
defines a ray $\ga\in\frR^j$ (we shall write
$\ga(t)=\exp(it\xi)$). Conversely, since $\tau(\ga(t))$ is bounded
and symmetric semigroup in $V_\tau$ for each ray $\ga\in\frR^j$
and any $\tau\in\Sp A$, we have
$$\frac{d}{dt}\tau(\ga(t))\Big|_{t=0}=i\tau(\xi_\ga)\leq0.$$
Thus, {\rm(\ref{coneg})} and {\rm(\ref{conex})} define the same
cone. If $\xi\in C^j$ and $-\xi\in C^j$, then $f(\exp(z\xi)j)$ is
bounded in $\bbC$ for each $f\in A_\fin$. Since this function is
continuous on $\bbC$ and analytic on $\bbC\setminus\bbR$, it is
entire. By Liouville Theorem, $f(\exp(z\xi)j)=f(j)$, but this
implies $\exp(z\xi)j=j$ since $A_\fin$ separates points of
$\cM_A$; hence, $\xi=0$. This proves that $C^j$ is pointed. The
remaining assertion is clear.
\end{proof}
We say that a ray $\ga$ is {\it $j$-central} if $\ga(0)=j$ and
$\ga(t)\in Z(G^j)$ for all $t\geq0$ (omitting $j$ in the notation
if $j=e$).
\begin{corollary}\label{cente}
If $C^j\neq0$, then the centre of $G^j$ is not discrete, and there
exists a nontrivial $j$-central ray $\ga$.
\end{corollary}
\begin{proof}
By Theorem~\ref{Liera}, $C^j\neq\{0\}$. If $\xi\in
C^j_A\setminus\{0\}$, then
$$\ze=\int\Ad(h)\xi\,dm_j(h)\neq0;$$
moreover, $\ze$ is $\Ad(G^j)$-fixed. Hence, $\exp(t\ze)$ lies in
the centre of $G^j$ and the semigroup $R_{\exp(it\ze)}R_j$
commutes with $R_h$ for all $h\in G^j$.
\end{proof}
By $G_{\bbC}$ we denote the complexification of the group $G$.
As a set, it can be identified with $G\exp(i\frg)$ (we may assume that
$G$ is a matrix group). Each finite
dimensional representation of $G$ uniquely extends to $G_{\bbC}$.
For $j\in\cI_A$ let $\frl^j$ be the real linear span of $C^j$ in $\frg^j$,
\begin{eqnarray*}
P^j=G^j\exp(iC^j).
\end{eqnarray*}
Since $C^j$ is $\Ad(G^j)$-invariant, $\frl^j$ is an ideal in $\frg^j$ and
\begin{eqnarray}\label{comgc}
G^j\exp(iC^j)=\exp(iC^j)G^j.
\end{eqnarray}
Let $L^j$, $L_{\bbC}^j$ be connected subgroups of $G^j$,
$G_{\bbC}^j$ that correspond to $\frl^j$, $\frl_{\bbC}^j$,
respectively. In the theorem below, they are equipped with the
underlying topology of a Lie group.

\begin{theorem}[\rm Analytic structure]\label{anast}
For each $j\in\cI_A$, the set $P^j$ is a subsemigroup of
$G^j_{\bbC}$, which consists of $g\in G^j_{\bbC}$ such that
\begin{eqnarray}\label{contr}
\|\tau(g)\|_{\BL(V^j_\tau)}\leq1
\end{eqnarray}
for all $\tau\in\Sp A$. The interior $D^j$ of
$L^j_{\bbC}\mathbin\cap P^j$ in $L^j_{\bbC}$ is an open semigroup
in $L^j_{\bbC}$, and the natural embedding $D^j\to\cM_A$ defines
the analytic structure, which is nontrivial if $C^j\neq0$.
\end{theorem}
\begin{proof}
Let $g=h\exp(i\xi)\in G^j_{\bbC}$, where $h\in G^j$,
$\xi\in\frg^j$. Then $\tau(\exp(i\xi))=\tau(g)\tau(h)^{-1}$.
Hence, $\tau(g)\in\tau(\cM_A)$ is equivalent to
$\tau(\exp(i\xi))\in\tau(\cM_A)$. Since $\tau(h)$ is isometric in
$V_\tau^j$, {\rm(\ref{contr})} holds for $\exp(i\xi)$ and for $g$
simultaneously. Further, {\rm(\ref{contr})} is true for
$\exp(i\xi)$ if and only if
\begin{eqnarray*}
\|\tau(\exp(it\xi))\|_{\BL(V^j_\tau)}\leq1
\quad\text{for all}\ t>0
\end{eqnarray*}
(the norm of a symmetric nonnegative operator is equal to its greatest
eigenvalue), and this 
is equivalent to the inequality in {\rm(\ref{conex})}. By
Theorem~\ref{Liera},  $P^j$ is distinguished by inequalities
{\rm(\ref{contr})}; hence, $P^j$ is a semigroup in $G^j_{\bbC}$.
Evidently, $P^j\mathbin\cap L^j_{\bbC}=L^j\exp(iC^j)$, and this
set has nonempty interior $D^j$ in $L^j_{\bbC}$, which is also a
semigroup. The embedding $D^j\to\cM_A$ defines the analytic
structure in $\cM_A$ because the mapping $g\to\tau(g)$ is
holomorphic on $L^j_{\bbC}$ for all $\tau\in\Sp A$, and the
uniform closure keeps this property. If $C^j\neq0$, then
$\frl^j\neq0$. Hence, $D^j$ is not a single point.
\end{proof}

\begin{lemma}\label{freei}
There exists an open subsemigroup $U$ of $\cM_A$ such that
$U\supseteq G$ and $(\vf^*\vf)^0=e$ for all $\vf\in U$.
\end{lemma}
\begin{proof}
Since $A_\fin$ is dense in $A$ and $A$ is separating, there exists
a finite dimensional representation $\tau$ of $G$ such that
$\Sp\tau\subseteq\Sp A$ and $\tau$ is faithful on $G$. By
Corollary~\ref{reext}, $\tau$ can be extended to $\cM_A$. Then
$\tau(j)\neq1$ for each $j\neq e$ in $\cI_A$ (this is a
consequence of (\ref{exten}), Proposition~\ref{idems}, 4), and the
choice of $\tau$). Hence, $\det\tau(j)=0$, and
$U=\{\vf\in\cM_A:\,\det\tau(\vf)\neq0\}$ satisfies the lemma.
Indeed, $\vf(\vf^*\vf)^0=\vf$ for all $\vf\in\cM_A$; since
$\tau(\vf)$  is invertible for $\vf\in U$, the idempotent
$(\vf^*\vf)^0$ is also invertible. This implies $(\vf^*\vf)^0=e$.
\end{proof}
\begin{proposition}\label{cgace}
If $A\neq C(G)$ then $C^e\neq0$.
\end{proposition}
\begin{proof}
Due to Theorem~\ref{antis}, the least idempotent $\eps$ does not coincide
with $e$. By Proposition~\ref{resta}, it is sufficient to prove
the assertion for $G_\eps$. According to Theorem~\ref{antis}, we may assume
that $A$ is antisymmetric on $G$. Since each function which is constant
on $G$ is constant on $\cM_A$, this means that $A$ is antisymmetric
on $\cM_A$.
Then $\cM_A$ is connected since $A$ contains the characteristic function
of any closed and open subset of $\cM_A$ 
by Shilov's Idempotent Theorem. It follows that
$U\mathbin\cap\left(\cM_A\setminus G\right)\neq\emptyset$ for any
neighbourhood $U$ of $e$ in $\cM_A$ (otherwise, $G$ is open in
$\cM_A$ due to the homogeneity).  Further, if $\vf\notin G$ then
$\vf^*\vf\notin G$ by Theorem~\ref{semit}. Thus, Lemma~\ref{freei}
and Theorem~\ref{rays.} imply the existence of a nontrivial ray
$\ga$ starting at $e$, and  Theorem~\ref{Liera} yields $C^e\neq0$.
\end{proof}
\begin{corollary}\label{ceray}
If $A\neq C(G)$, then there exists $\ze\in\Int C^e$ such that the
corresponding one parameter group in $G$ is a subgroup of $Z(G)$
and $\exp(it\ze)$, $t\geq0$, is a nontrivial ray in $Z(\cM_A)$.
\end{corollary}
\begin{proof}
Combine Proposition~\ref{cgace}, Corollary~\ref{cente}
and Proposition~\ref{centr}.
\end{proof}
\begin{lemma}\label{jkdim}
If $j,k\in\cI_A$, $j\leq k$ and $j\neq k$, then $\dim G_j>\dim
G_k$.
\end{lemma}
\begin{proof}
If $\dim G_j=\dim G_k$, then the homogeneous space $G_j/G_k$ is
finite and the algebra $R_kA\Big|_{G_j}$ is separating on it by
Proposition~\ref{averr} and Proposition~\ref{idems}, 3). Each
separating algebra on a finite set is evidently the algebra of all
functions on it.  Hence, the measure $m_j$ cannot be
multiplicative on $A$, contradictory to Proposition~\ref{idems},
4).
\end{proof}
\begin{lemma}\label{j-ray}
Let $j\in\cI_A$. If $j\neq\eps$, then $Z(G^j)$ contains a
nontrivial ray starting at $j$. Moreover, if $j\in Z(\cM_A)$ and
$j\neq\eps$, then there exists a nontrivial ray lying in
$Z(\cM_A)$ and starting at $j$.
\end{lemma}
\begin{proof}
If $\dim G=0$, then $G$ is finite. Since $A$ is separating,
$A=C(G)$; hence, $\cM_A=G$ and $e=\eps$. This provides the base
for the induction by $\dim G$. Let $\dim G>0$, $e\neq\eps$, and
$\ga$ be a nontrivial central ray in $\cM_A$ starting at $e$
(existing by Corollary~\ref{ceray}). Then $k=\ga(\infty)\neq e$.
If $k\not\geq j$,  then $\ga(t)j$ is a ray in $Z(G^j)$; moreover,
it is contained in $Z(\cM_A)$ if $j\in Z(\cM_A)$. This ray is
nontrivial because $kj\neq j$. Let $k\geq j$. Since $\ga$ is
central, $k\in Z(\cM_A)$. Further, $\dim G_k>0$ according to
Lemma~\ref{jkdim}. Therefore, $\dim G/G_k<\dim G$ and the
induction hypothesis can be applied to the algebra $R_kA$ on
$G^k=G/G_k$, whose maximal ideal space can be identified with the
semigroup $k\cM_A$ by Theorem~\ref{surje}. It follows that either
$j=\eps$ or $\frR^j$ contains a nontrivial ray $\td\ga$ in
$Z(G^j)\mathbin\cap k\cM_A$. If $j\in Z(\cM_A)$, then $\td\ga$ can
be assumed to be central in $k\cM_A$ by Corollary~\ref{ceray};
then it is evidently central in $\cM_A$.
\end{proof}
\begin{corollary}
$C^j=0$ if and only if $j=\eps$.\qed
\end{corollary}
We say that a finite sequence $\ga_1,\dots,\ga_n$ is {\it a chain
of rays of the length $n$} if all rays are nontrivial and
\begin{eqnarray*}
\ga_l(\infty)=\ga_{l+1}(0),\quad l=1,\dots,n-1.
\end{eqnarray*}
Also, we say that the chain joins $\ga_1(0)$ and $\ga_n(\infty)$.
\begin{theorem}\label{seque}
For any pair $j,k\in\cI_A$ such that $j\leq k$, there exists a
chain of rays which joins them. Its length does not exceed $\dim
G_j-\dim G_k$.
\end{theorem}
\begin{proof}
It is sufficient to prove the proposition for the restriction of
$A$ to $G_j$. Hence, we may assume $j=\eps$ but in this case the
assertion is almost evident: by Lemma~\ref{j-ray}, there exists a
chain of rays starting at $k$, and it cannot be continued only if
its endpoint is the least idempotent $\eps$. The upper bound of
the length follows from Lemma~\ref{jkdim}.
\end{proof}
Let $\si_j$ denote the Haar measure of $Z_j$.
\begin{proposition}\label{recen}
For any $j\in\cI_A$, $\si_j\in\cM_j$.
\end{proposition}
\begin{proof}
It is sufficient to prove that $\si_\eps\in\cM_\eps$,  assuming
$G=G_\eps$ (then we get the proposition applying this to
$A\big|_{G_j}$). We use the induction on $\dim G$.  The assertion
is evident if $\dim G=0$. If $\dim G>0$, then there exists a
nontrivial central ray $\ga\in\frR^e$. For all $f\in A$ and $t>0$
\begin{eqnarray*}
f(\ga(t))=\frac1\pi\int_{-\infty}^\infty\frac{tf(\ga(ix))}{x^2+t^2}\,dx
\end{eqnarray*}
due to Theorem~\ref{Liera}. Hence, $\cM_{\ga(t)}$ contains a
measure concentrated on a central one-parameter subgroup of $G$.
Therefore, $k=\ga(\infty)$ has a representing measure with the
support in its closure that is a connected central subgroup of
$G$. As in Theorem~\ref{seque}, $k\in Z(\cM_A)$,
$\cM_{R_kA}=k\cM_A$, $G^k=G/G_k$ and $\dim G^k<\dim G$. By the
induction hypothesis, the Haar measure $\td\si_\eps$ of the
identity component $\td Z_\eps$ of $Z_{G^k_\eps}(G^k_\eps)$ is
representing for $\eps$ on $G^k$. We claim that $\td
Z_\eps=kZ_\eps$. Indeed, $\td Z_\eps=\exp(\td{\frz}_\eps)$ is the
exponent of the set of $\Ad(G^k)$-fixed points in $\frg^k$;
$\pi_k$-preimages in $\frg$ of all these points are
$\Ad(G)$-invariant, hence contain $\Ad(G)$-fixed points.
Therefore, $\td\si_\eps=\pi_k\si_\eps$. This implies
$\si_\eps\in\cM_\eps$.
\end{proof}
\begin{proposition}\label{absup}
For any $s\in\cS_A$, there exists a connected abelian subgroup
$H\subseteq G$ such that
\begin{eqnarray*}
s\in\wh H\subseteq Z(s).
\end{eqnarray*}
\end{proposition}
\begin{proof}
Let $j=s^0$, $\nu$ be the Poisson measure
$\frac1\pi\frac{dt}{t^2+1}$ on $\bbR$. By Proposition~\ref{recen},
Theorem~\ref{rays.} and {\rm(\ref{onepa})}, for some
$\xi\in\frn_j$
\begin{eqnarray*}
f(s)=\int_{\bbR}f(s^{it})\,d\nu(t)=
\int_{\bbR}f(j\exp(t\xi))\,d\nu(t)=\phantom{xxxxxxxxxxxxxxxxxxxx}\\
\int_{\bbR}\int_{Z_j}f(h\exp(t\xi))\,d\nu(t)\,d\si_j(h).
\end{eqnarray*}
Therefore, $s$ has a representing measure concentrated on the
closure $H$ of the set $Z_j\exp(\bbR\xi)$. Clearly, $Z_j$ is
$\Inn(N_j)$-invariant; hence, the inclusion $\xi\in\frn_j$ implies
that $H$ is a group. Furthermore, $H$ is abelian since it is
compact and solvable. Evidently, $H$ is connected. If $h\in
Z_j\subseteq G_j$, then $hs=sh=s$ due to equalities $s=js=sj$ and
$jh=hj=j$; if $h=\exp(t\xi)$, $t\in\bbR$, then $hs=s^{t+1}=sh$.
Thus, $H\subseteq Z(s)$.
\end{proof}
Proposition~\ref{absup} enables to prove a stronger version of
Theorem~\ref{antis} for connected Lie groups.
\begin{theorem}\label{antor}
Let $G$ be connected, $T$ be a maximal torus in $G$, $\tau$ be the
Haar measure of{\,} {$T$}, and $A$ be an invariant algebra on $G$.
Then the following assertions are equivalent:
\begin{itemize}
\item[\rm a)] $A$ is antisymmetric;
\item[\rm b)] the closure of $A\big|_T$ in $C(T)$ is antisymmetric;
\item[\rm c)] $\tau$ is multiplicative on $A$.
\end{itemize}
In this case $\tau\in\cM_\eps$, where $\eps$ is the least idempotent.
\end{theorem}
\begin{proof}
Let $A$ be antisymmetric. Then $G=G_\eps$, and $Z_\eps$ is the
identity component of the centre of $G$. Hence, $Z_\eps\subseteq
T$. By Proposition~\ref{recen}, $\si_\eps\in\cM_\eps$. Due to
Theorem~\ref{antis}, $\eps$ is zero of $\cM_A$; in particular,
$\de_t*\si_\eps\in\cM_\eps$ for all $t\in T$. The convex weakly
closed span of these measures contains $\tau$. This proves the
implication a)$\ \Rightarrow\ $c) and the last assertion of the
theorem. If $A$ is not antisymmetric, then it contains a
nonconstant real function $f$. Let $f(g)\neq f(h)$, where $g,h\in
G$. Since $G$ is connected, there exist $x,y\in G$ such that
$xgy,xhy\in T$. Then $L_x^{-1}R_y^{-1}f\in A$ is a real function,
which is not constant on $T$. Thus, b) implies a). It remains to
note that b) and c) are equivalent by Theorem~\ref{antis}.
\end{proof}
For Lie groups, the ordered set $\cI_A$ has some special properties.
\begin{proposition}
For any $k\in\cI_A$, there exists its neighbourhood $U$ in $\cM_A$
such that $j\in\cI_A\mathbin\cap Z(k)\mathbin\cap U$ implies
$k\leq j$.
\end{proposition}
\begin{proof}
If $jk=kj$, then $jk\in\cI_A$. Let $\rho$ be a finite dimensional
representation of $G$ such that $\Sp\rho\subset\Sp A$ and $\rho$
is faithful on $G^k$. Set
\begin{eqnarray*}
W=\left\{\vf\in\cM_A:\,\rank\rho(\vf)\geq \rank\rho(k)\right\}.
\end{eqnarray*}
Then the interior of $W$ contains $k$. For any nontrivial ray
$\ga\in\frR^k$, the one-parameter semigroup $\rho(\ga(t))$ is
nontrivial due to Theorem~\ref{Liera}. Hence, $\ga(\infty))\notin
W$. Taken together with Theorem~\ref{seque}, this implies that $W$
does not contain $j\in\cI_A$ if $j\leq k$ and $j\neq k$. Let
$j\in\cI_A\mathbin\cap Z(k)$. Then $jk\in\cI_A$ and $jk\leq k$.
Therefore, $jk\neq k$ implies $jk\notin W$. Thus, any
neighbourhood $U\subseteq W$ of $k$ such that $U^2\subseteq W$
satisfies the proposition.
\end{proof}
\begin{proposition}\label{Lieid}
Let $j\in\cI_A$, and let $T$ be a maximal torus in $G$. Each of
the following conditions is equivalent to the inclusion $j\in\wh
T$:
\begin{itemize}
\item[a)] $j\in Z(T)$;
\item[b)] $T\subseteq N_j$;
\item[c)] $T\mathbin\cap G_j$ is a maximal torus in $G_j$;
\item[d)] $Z_j\subseteq T$.
\end{itemize}
In this case, $Tj$ is a maximal torus in $G^j$ and
\begin{eqnarray}
&Gj\mathbin\cap Z(T)\subseteq G^j.\label{gjzt.}
\end{eqnarray}
\end{proposition}
\begin{proof}
According to Proposition~\ref{resta}, $\wh T$ is abelian. Hence,
$j\in\wh T$ implies a). By a), $T\subseteq Z(j)=N_j$, i.e., b)
holds. If $T\subseteq N_j$, then $H=TG_j$ is a subgroup of $G$,
which is locally a direct sum of a torus $T'\subseteq T$ and
$G_j$. Then $T$ and $T/T'$ are maximal tori in $H$ and $H/T'$,
respectively. Since $H/T'$ is locally isomorphic to $G_j$ and all
maximal tori are conjugated, the torus $T\mathbin\cap G_j$ is
maximal in $G_j$. Thus, c) follows from b). The implication c)$\
\Rightarrow\ $d) is trivial. If d) is true, then $j\in\wh T$ by
Proposition~\ref{recen} and Lemma~\ref{rpeak}.

Clearly, $T$ is a maximal torus in $N_j$. Thus, b) implies that
$Tj\cong T/(T\mathbin\cap G_j)$ is a maximal torus in $G^j\cong
N_j/G_j$.

If $t\in T$, $g\in G$ and $gj\in Z(T)$, then $tgj=gjt$. By d),
\begin{eqnarray*}
f(jgj)=\int f(tgj)\,d\si_j(t)=\int f(gjt)\,d\si_j(t)=f(gj),\quad f\in A,
\end{eqnarray*}
where $\si_j$ is the measure of Proposition~\ref{recen}.
This implies $jgj=gj$ and $g^{-1}jgj=j$.
Applying Lemma~\ref{invma} and Proposition~\ref{idems},
3) and 4), we get $g^{-1}G_jgG_j=\supp(\de_{g^{-1}}*m_j*\de_g*m_j)=G_j$.
Therefore, $g^{-1}G_jg\subseteq G_j$. i.e. $g\in N_j$ and $gj\in G^j$.
This proves {\rm(\ref{gjzt.})}.
\end{proof}

\begin{proposition}
Let $T$ be a maximal torus in $G$. Then
\begin{eqnarray}
\cS_A\mathbin\cap Z(T)=\cS_A\mathbin\cap\wh T. \label{cziat}
\end{eqnarray}
\end{proposition}
\begin{proof}
Let $s\in\cS_A\mathbin\cap Z(T)$. Then $s^{z}\in Z(T)$ for all
$z\in\bbC^+$ due to Corollary~\ref{comst}. In particular,
$j=s^0\in Z(T)$. By Proposition~{\rm\ref{Lieid}}, a),
$j\in\wh T$. Thus, one parameter group $\ga(t)=s^{it}$, $t\in\bbR$,
lies in $G^j\mathbin\cap Z(T)\subseteq Z(Tj)$.
Since $Tj$ is a maximal torus in $G^j$ by Proposition~\ref{Lieid},
$\ga(\bbR)\subseteq Tj$. It follows from Theorem~\ref{Liera} that
$s\in\wh{\ga(\bbR)}$; hence $s\in\wh{Tj}\subseteq\wh T$.
This proves the nontrivial inclusion
$\cS_A\mathbin\cap Z(T)\subseteq\cS_A\mathbin\cap\wh T$.
\end{proof}
The following theorem reduces many problems concerning
invariant algebras to the abelian case. We shall say that
the semigroup $\wh T$ below is a {\it Cartan subsemigroup of $\cM_A$}.
\begin{theorem}\label{Carta}
Let $T$ be a maximal torus in $G$.
\begin{itemize}
\item[\rm a)] For each $s\in\cS_A$, there exists $g\in G$ such
that $g^{-1}sg\in\wh T$. \item[\rm b)] Let $\vf\in Z(\wh T)$ and
let $\vf=gs$, where $s\in\cS_A$, $g\in G$, be its polar
decomposition. Then $s\in\wh T$ and $gj\in G^j\mathbin\cap Z(T)$,
where $j=s^0\in\wh T$. \item[\rm c)] If $G^j$ is connected for
every $j\in\cI_A$, then $\wh T$ is a maximal abelian subsemigroup
of $\cM_A$.
\end{itemize}
\end{theorem}
\begin{proof}
Since any connected abelian subgroup of $G$ is conjugated with a subgroup
of $T$ by an inner automorphism, a) follows from
Proposition~\ref{absup} and Lemma~\ref{rpeak}.

If $\vf\in Z(\wh T)$, then $\vf^*\in Z(\wh T^*)=Z(\wh T)$, hence,
$\vf^*\vf\in Z(\wh T)$. By Corollary~\ref{comst}, $(\vf^*\vf)^t\in
Z(T)$ for all $t\geq0$. In particular, $s=(\vf^*\vf)^{\frac12}$
and $j=s^0$ belong to $Z(T)$. Applying {\rm(\ref{cziat})}, we get
$s,j\in\wh T$. Since $R_s$ and $R_j$ are self-adjoint in $H^2$,
equalities $sj=js=s$ imply $\ker R_s=\ker R_j$. Thus, the operator
$R_s^{-1}$ is well-defined on $R_jA_\fin$. Therefore, for all
$f\in R_jA_\fin$ and $t\in T$
\[
R_t R_{gj}f=R_tR_\vf R_s^{-1}f=R_\vf R_s^{-1}R_tf=R_{gj}R_tf.
\]
Hence, $gj\in Z(T)=Z(\wh T)$, and a reference to
{\rm(\ref{gjzt.})} concludes the proof of b).

If $G^j$ is connected, then $Tj$ is a maximal abelian subgroup of
$G^j$ being its maximal torus. Therefore, $Z(\wh T)\mathbin\cap
G^j=Tj\subseteq\wh T$. Thus, b) and the assumption of c) imply
$Z(\wh T)=\wh T$.
\end{proof}
Of course, $\wh T$ need not be a maximal abelian semigroup in
$\cM_A$ in general. For instance, this is true if $G$ is abelian
disconnected. An illuminating example is the group
$G=\bbT\times\bbZ_2$ and the invariant algebra $A$ corresponding
to the semigroup
\[
S=\{(n,1):\,n\in\bbZ,\,n\geq0\}\mathbin\cup\{(n,-1):\,n\in\bbZ,\,n>0\}.
\]
in the group $\wh G=\bbZ\times\bbZ_2$, where $\bbZ_2$ is realized
as the multiplicative group $\{\pm1\}$. Also, $A$ can be described
as the algebra of pairs of functions $f_1,f_2\in A(\bbD)$ (see
Example~\ref{disca}) such that $f_1(0)=f_2(0)$; $\cM_A$ is the
union of two discs, which have a single common point.

\section{Invariant algebras on tori}
Theorem~\ref{abelc} identifies $\cM_A$ and $\Hom(\Sp A,\bbD)$. If
$G=\bbT^n$, then $\wh G=\bbZ^n$, $\Sp A$ is a semigroup in
$\bbZ^n$, and the algebraic structure of $\Hom(\Sp A,\bbD)$ is
determined by the geometry of $\Sp A$. Main ingredients are the
asymptotic cone $\al(\Sp A)$ (see {\rm(\ref{acone})}) and the face
structure. Everywhere in this section, $\Sp A$ is denoted by $S$;
it can be an arbitrary semigroup in $\bbZ^n$. We start with a
description of subsemigroups of $S$ corresponding to idempotents
in $\Hom(S,\bbD)$. They will be called {\it faces of $S$}. We
shall prove that one can get the collection $\frI_S$ of all faces
by a procedure of step-by-step expansion. It starts with
$\frI_S=\{S\}$ and finishes when $\frI_S$ does not change. On each
step, the following operation is applied to the faces $P$ that are
already found:
\begin{eqnarray}\label{proce}
\text{if}\   F\in\frF_{\al(P)}\  \text{and}\
P\mathbin\cap F\neq\emptyset\  \text{then add}\
P\mathbin\cap F\  \text{to}\   \frI_S.
\end{eqnarray}
Obviously, the procedure is finite.
Since $\cM_A$ is abelian, the set $\cS_A$ is a semigroup.
According to {\rm(\ref{invab})},
it can be identified with $\Hom(S,\bbI)$, where $\bbI=[0,1]$
is considered as a multiplicative semigroup.
Each idempotent in $\Hom(S,\bbI)$ corresponds to the characteristic
function $\ka_X$ of some set $X\subseteq S$:
\begin{eqnarray}\label{charf}
\ka_X(x)=\left\{\begin{array}{lr}1,&x\in X,\cr
                               0,&x\notin X.\end{array}\right.
\end{eqnarray}
We denote $\ka_S$ by ${\bf 1}$.
\begin{lemma}\label{rayho}
Let $\ga$ be one parameter semigroup in $\Hom(S,\bbI)$. Suppose that
$\ga(0)={\bf 1}$. Then $\ga(t)(x)=e^{-t\la(x)}$, where $\la\in S^\star$,
$t\geq0$, $x\in S$. Furthermore, $\ga(\infty)=\ka_P$,
where $P=S\mathbin\cap F$ for some $F\in\frF_{\al(S)}$.
\end{lemma}
\begin{proof}
Without loss of generality we may assume that $S$ linearly
generates $\bbR^n$. Let $x\in S$. The assumption $\ga(0)(x)=1$,
taken together with the continuity and the semigroup property,
implies $\ga(t)(x)>0$ for all $t\geq0$. Hence, a real valued
function
$$\la_t(x)=-\log\ga(t)(x)$$
is well defined. Clearly, $\la_t$ is additive on $S$. Thus, it
extends to the additive functional on the cocompact discrete group
$S-S$, hence to the linear functional on $\bbR^n$. Further,
$\la_t=t\la_1$ since $\la_0(x)=0$ and $\la_t(x)$ is an additive
and continuous function on $t$ for any $x\in S$. Obviously,
$\la_t\in S^\star=\al(S)^\star$. Therefore, $\ga(t)=e^{-t\la}$,
where $\la=\la_1\in S^\star$. It follows that $\ga(\infty)(x)=0$
if $\la(x)>0$, and $\ga(\infty)(x)=1$ if $\la(x)=0$; in other
words, $\ga(\infty)=\ka_P$, where
$$P=S\mathbin\cap F,\qquad F=\al(S)\mathbin\cap\la^{-1}(0).$$
It remains to note that $F\in\frF_{\al(S)}$ since $\la\in
S^\star$.
\end{proof}
\begin{corollary}\label{opens}
The set $\cS_A^e$ is open in $\cS_A$; moreover, for any
$x\in\Int\al(S)$
$$\cS_A^e=\{\chi\in\Hom(S,\bbI):\,\chi(x)\neq0\}.$$
\end{corollary}
\begin{proof}
Let $\chi\in\Hom(S,\bbI)=\cS_A$. By definition, if
$\chi\notin\cS_A^e$, then $\chi\ka_P=\chi$ for some $P\neq S$. By
Lemma~\ref{rayho} and Theorem~\ref{seque}, this implies $\chi=0$
on $\Int\al(S)\mathbin\cap S$.
\end{proof}
\begin{proposition}\label{sefac}
For a set $P\subseteq S$, the inclusion $P\in\frI_S$ is equivalent
to each of following conditions:
\begin{itemize}
\item[\rm 1)] $P$ is a semigroup and $S\setminus P$ is an ideal
of the semigroup $S$;
\item[\rm 2)] $\ka_P\in\Hom(S,\bbD)$.
\end{itemize}
Moreover, if $P\in\frI_S$ and $P\neq S$, then $\al(P)$ is a closed
cone containing in some face $F\in\frF_{\al(S)}$, $F\neq\al(S)$.
\end{proposition}
\begin{proof}
Let $P\in\frI_S$. Proving 1), we may assume without loss of
generality that $P$ is added to $\frI_S$ on the first step, i.e.,
that $P=F\mathbin\cap S$ for some closed face $F$ of the cone
$\al(S)$. Then $P$ is a semigroup and $S\setminus P$ is an ideal
since $F$ has these properties with respect to $\al(S)$. The
equivalence of 1) and 2) is obvious. Taken together with
Proposition~\ref{sefac} and Lemma~\ref{rayho}, 2) implies
$P\in\frI_S$. The last assertion follows from {\rm(\ref{proce})}.
\end{proof}
The semigroup $S$ uniquely determines $\cM_A$ but the converse is
not true (for example, the semigroup in $\bbZ$ which is generated
by numbers $0, 2, 3$ corresponds to the same maximal ideal space
as the semigroup of all nonnegative integers).
Proposition~\ref{sefac} and Lemma~\ref{rayho} show how to
reconstruct $\cM_A$ by $S$. We outline this procedure (without
proofs) and construct a kind of a hull for $S$ which corresponds
to the same maximal ideal space as $S$. Besides, we give some
examples to illustrate the complexity of the object.

Topologically, $\cM_A$ is $G\times\cS_A$ factorized by some
equivalence. By Theorem~\ref{polar}, to decide if $(g,s)$ and
$(h,s)$ represent the same point in $\cM_A$ or not one have to
know if $g^{-1}h\in G^j$, where $j=s^0$; the component $s$ must
coincide. According to Theorem~\ref{Liera}, $\cS_A$ is the union
of cones $C^j$, $j\in\cI_A$, which are compactified and glued
together. Thus, to build $\cM_A$, one needs
\begin{eqnarray}\label{datma}
\cI_A,\quad\frC_A=\{C^j:\,j\in\cI_A\},\quad\frG_A=\{G^j:\,j\in\cI_A\},
\end{eqnarray}
and has to know points which must be identified. The geometric
object which contains information on $\cS_A$ (precisely, $\cI_A$,
$\frC_A$, and the topology in the union of $C^j$, $j\in\cI_A$) is
the family of asymptotic cones to faces of $S$:
\begin{eqnarray*}
\frC_S=\left\{\al(P):\,P\in\frI_S\right\}.
\end{eqnarray*}
Indeed, faces are in one-to-one correspondence with idempotents in
$\cM_A$ by Proposition~\ref{sefac}, 2) and Theorem~\ref{abelc}.
For each $j\in\cI_A$, the cone $C^j\subseteq\frg^j$ of
Theorem~\ref{Liera} is dual to $C_j\in\frC_S$. The topology in the
union of these cones relates to the topology of the pointwise
convergence in $\Hom(S,\bbI)$. It is sufficient to describe the
closure of $C^e=C_e^\star$, where $C_e=\al(S)$ (the operation must
be repeated for all $j\in\cI_A$). The first step of the procedure
{\rm(\ref{proce})} defines idempotents lying in the closure of
$\cS_A^e$ in $\cM_A$. The set $\cS_A^e$, by Theorem~\ref{Liera},
may be identified with $C^e$ and, by Lemma~\ref{rayho}, with the
set
$$\{e^{-\la}:\,\la\in C^e\}\subseteq\Hom(S,\bbI).$$
Let $j\in\cI_A$ belong to the closure of $\cS_A^e$ and $P$ be the
corresponding face of $S$. Clearly, $C^j$ contains all functionals
in $C_e^\star$ restricted to $C_j$ but it can be wider if
$$C_j=\al(P)$$
is a proper subset of the closed face of $C^e$ that
includes $P$. The restriction of $C_e^\star$ can be identified
with the projection $\pi_j(C^e)$; the closure of $C^e$  is the
union of these sets. Continuing this, we get $\cS_A$.

To find $\frG_A$, it is sufficient to know for each $j\in\cI_A$
the annihilator $\Ga_j$ of $G^j$ consisting of $y\in\wh G$
such that $y(g)=1$ for all $g\in G_j$. Then
\begin{eqnarray}\label{anngj}
G^j=\{g\in G:\,y(g)=1\ \text{for all}\  y\in\Ga_j\}
\end{eqnarray}
due to Pontrjagin duality. Clearly, $\Ga_j$ may be replaced by
$C_j\mathbin\cap S$ in {\rm(\ref{anngj})} and $\Ga_j$ is the group
generated by $C_j\mathbin\cap S$.

Let $\frG_S$ be the collection of groups $\Ga_j$, $j\in\frI_S$.
Families
\begin{eqnarray}\label{sdata}
\frI_S,\ \frC_S,\ \frG_S
\end{eqnarray}
correspond to $\cI_A,\frC_A,\frG_A$, respectively, and uniquely
determines them. Objects {\rm(\ref{sdata})} satisfy following conditions:
\begin{itemize}
\item[A)] if $j\leq k$ then $C_j$ is a closed cone which is contained
in some face of $C_k$;
\item[B)] $C_j=\al(\Ga_j\mathbin\cap C_j)$;
\item[C)] if $j\leq k$ then $\Ga_j\subseteq\Ga_k$,
\end{itemize}
where $j,k\in\cI_A\cong\frI_S$. Set
\begin{eqnarray}\label{envel}
\wtd S=\bigcup\nolimits_{j\in\frI_S}\wtd P_j, \quad\text{where}\quad
\wtd P_j=
\Ga_j\mathbin\cap\Int C_j.
\end{eqnarray}
Then $\wtd S\supseteq S$ and $\wtd{\!\wtd S}=\wtd S$. The set
$\wtd S$ is a semigroup, and families {\rm(\ref{sdata})} for $\wtd
S$ are the same as for $S$. Thus, $\cM_{\wtd A}=\cM_A$, where
$\wtd A$ is the invariant algebra corresponding to $\wtd S$.
Combined with the procedure {\rm(\ref{proce})}, conditions A)--C)
give an universal method to construct maximal ideal spaces of
invariant algebras on tori.

The case of finitely generated semigroups is simpler than the
general one. We keep the notation of Theorem~\ref{anast};
``compactification" means ``homeomorphic embedding to a compact
space as a dense open subset".
\begin{proposition}\label{fings}
Suppose that the semigroup $S\subseteq\bbZ^n$ is generated by
its finite subset. Then
\begin{itemize}
\item[\rm a)] the correspondence between $\frI_S$ and
$\frF_{\al(S)}$ defined by $P=F\mathbin\cap\al(S)$, where
$P\in\frI_S$ and $F\in\frF_{\al(S)}$, is one-to-one; \item[\rm b)]
$\cS_A^e$ is dense and open in $\cS_A$; \item[\rm c)] if $A$ is
antisymmetric and separating, then $G=L^e$, and $\cM_A$ is a
compactification of the domain $D^e$ in the complex torus
$G_{\bbC}$.
\end{itemize}
\end{proposition}
\begin{proof}
Set $C=\al(S)$ and let $X=\{x_1,\dots,x_m\}$ be the finite set generating
$S$. Then
\begin{eqnarray}\label{finco}
S=\bbN x_1+\dots+\bbN x_m,\quad
C=\bbR^+x_1+\dots+\bbR^+x_m.
\end{eqnarray}
Each closed face $F\in\frF_C$ of the convex closed cone $C$ has
the form $F=\la^\bot\mathbin\cap C$ for some $\la\in C^\star$
(precisely, for every $\la\in\Int(F^\bot\mathbin\cap C^\star)$).
Set $X_\la=\la^\bot\mathbin\cap X$. Then
\begin{eqnarray}\label{facef}
F=\sum_{x\in X_\la}\bbR^+x
\end{eqnarray}
since $\la(x)>0$ for $x\notin X_\la$. This implies that $F=\al(P)$
for $P=F\mathbin\cap S$ and that
\begin{eqnarray}\label{frfff}
\frF_F\subseteq\frF_C.
\end{eqnarray}
Therefore, the mapping $P\to\al(P)$ is inverse to $F\to
F\mathbin\cap S$ and a) is true. It follows from
{\rm(\ref{facef})} that the projection $C^\star\to F^\star$, dual
to the embedding $F\to C$, is surjective. Combining this with
Lemma~\ref{rayho}, we get that each one parameter semigroup $\ga$
in $\Hom(P,\bbI)$ such that $\ga(0)=\ka_P$ can be realized in the
form $\ka_P\td\ga(t)$, where $\td\ga$ is one parameter semigroup
in $\Hom(S,\bbI)$.  Thus, {\rm(\ref{frfff})}, Theorem~\ref{rays.},
Theorem~\ref{seque}, and Corollary~\ref{opens} yields b). Due to
Theorem~\ref{polar} and b), the set $P^e=G\cS_A^e$ is dense and
open in $\cM_A$ (in fact, this is the set $U$ of
Lemma~\ref{freei}). If $A$ is antisymmetric, then $C$ is pointed
according to {\rm(\ref{antab})} and {\rm(\ref{finco})}. Hence, the
dual cone $C^\star=C^e\subset\frg$ has nonempty interior in
$\frg$; consequently $\frl=\frg$.  Since $G$ is connected and $A$
is separating, $G=L^e$ and $D^e$ is the interior of $P^e$ in the
group $G_{\bbC}$. The polar decomposition on the set $P^e$ is
one-to-one, hence locally homeomorphic. Therefore, $D^e$ is dense
in $P^e$; this proves c).
\end{proof}
\begin{corollary}\label{finid}
If $S$ is finitely generated, then the set $\cI_A$ is finite.
\end{corollary}

\section{Examples}
We shall say that an invariant algebra $A$ is {\it finitely
generated} if some its finite dimensional invariant subspace
generates it as a Banach algebra. If $G$ is abelian, then this is
equivalent to the condition that the semigroup $S=\Sp A$ is
finitely generated. It is not difficult to prove that $A$ is
finitely generated if and only if there exist a finite dimensional
Hilbert space $V$ and an isomorphic embedding of $G$ to $\rU(V)$
such that $A_\fin=\cP(\rU(V))\big|_G$; then $\cM_A$ coincides with
the polynomially convex hull $\wh G$ of $G$ in $\BL(V)$. Further,
if $A$ is finitely generated, then the closure of the restriction
$A\big|_T$ to the maximal torus $T$ is a finitely generated
invariant algebra on $T$ (probably, the converse is also true). If
$G$ is abelian, then the compactification of $D^e$ in
Proposition~\ref{fings} can be identified with its closure in the
affine variety $\Hom(S,\bbC)$, which admits a realization as an
algebraic submanifold of $\bbC^n$. Precisely,
\begin{eqnarray*}
\clos D^e=\{\chi\in\Hom(S,\bbC):\,|\chi(x)|\leq1\ \text{for all}\
x\in X\},
\end{eqnarray*}
where $X$ is a finite generating set. The example below
demonstrates some typical effects that occur in the realization
above.
\begin{example}\rm
Conditions $p,q\geq0$, where $p$ is even if $q=0$, distinguish a
semigroup $S\subset\bbZ^2$. It is generated by $x_1=(2,0)$,
$x_2=(0,1)$, $x_3=(1,1)$, which  satisfy the relation
$x_1+2x_2=2x_3$. The algebra $A$ is isomorphic to the closed (in
the $\sup$-norm on $\bbT^2$) linear span of all holomorphic
monomials $z_1^pz_2^q$ in $\bbC^2$ except odd powers of $z_1$.
There are $4$ idempotents in $\cI_A$ corresponding to $4$ faces of
$S$ (two of them are trivial). Clearly, $\frI_S=\frI_{\td S}$ and
$\frC_S=\frC_{\td S}$ for the semigroup $\td
S=\{(p,q)\in\bbZ^2:\,p,q\geq0\}$ but the group $G^j$ corresponding
to the face $p\geq0,q=0$, is disconnected; $\cM_A$ can be realized
as the bidisc $|z_1|,|z_2|\leq1$ in $\bbC^2$ with identified
points $(\pm z_1,0)$. Putting $z_k=\chi(x_k)$, $k=1,2,3$, where
$\chi\in\Hom(S,\bbD)$, we get an embedding of $\cM_A$ to the
variety $z_1z_2^2=z_3^2$  in $\bbC^3$; the image coincides with
that of the mapping $\bbD^2\to\bbC^3$,
$(\ze_1,\ze_2)\to(\ze_1^2,\ze_2,\ze_1\ze_2),$ and is distinguished
by inequalities $|z_1|,|z_2|,|z_3|\leq1$.\qed
\end{example}
\noindent Following three examples illuminate some properties of
infinitely generated algebras, which finitely generated ones do
not possess. According to the first of them, $\cS_A^e$ need not be
dense in $\cS_A$. This is a modification of an example studied in
\cite[Section~5]{Ar} with many details.
\begin{example}\label{notde}
\rm Let $S$ be the semigroup of $(p,q,r)\in\bbZ^3$ such that
either $r>0$ or $r=0$ and $p,q\geq0$. Then $\frC_S$ consists of
$5$ cones: $r\geq0$; $r=0$, $p,q\geq0$; $r=p=0$, $q\geq0$;
$r=q=0$, $p\geq0$; $r=p=q=0$. The set $\cS_A$ is homeomorphic to
the union of the unit disc $\bbD$ and the interval $[1,2]$ which
represents $\cS^e$ (with $2$ corresponding to $e$). Hence
$\cS_A^e$ is not dense in $\cS_A$ (cf. Proposition~\ref{fings},
b)). The semigroup $S$ is not finitely generated but each
character $\chi\in\Hom(S,\bbD)$ is uniquely determined by its
values at points $x_1=(1,0,0)$, $x_2=(0,1,0)$, $x_3=(0,0,1)$.
Indeed, for every $x\in S$ and all sufficiently large $n,m$ vector
$x+nx_1+mx_2$ belongs to the semigroup generated by $x_1,x_2,x_3$.
This uniquely defines $\chi(x)$ if $\chi(x_k)\neq0$, $k=1,2$. Set
$z_k=\chi(x_k)$, $k=1,2,3$. If $z_k=0$ for $k=1$ or $k=2$, then
$\chi(x+x_k)=0$ for any $x\in S$; hence, $\chi(p,q,r)=0$ for all
$p,q\in\bbZ$ and $r\geq1$ but $z_1,z_2$ uniquely determine
$\chi(p,q,0)$. Further, if $z_3\neq0$, then $|z_1|=|z_2|=1$:
otherwise, $\chi$ cannot be bounded on $\bbZ x_1+\bbZ x_2+ x_3$.
If $z_3=0$, then there is no restriction on $z_1,z_2$ except
$z_1,z_2\in\bbD$. Thus, $\cM_A=\cM_1\cup\cM_2$, where
\begin{eqnarray*}
&\cM_1=\{(z_1,z_2,z_3):\,|z_3|\leq1,\ |z_1|=|z_2|=1\},\\
&\cM_2=\{(z_1,z_2,z_3):\,z_3=0,\ |z_1|\leq1\}.
\end{eqnarray*}
Both sets have real dimension 4. Functions of $A$ are analytic on
$z_3$ in $\cM_1$ and on $z_1,z_2$ in $\cM_2$; $A$ is
antisymmetric.\qed
\end{example}
\noindent If $S$ is not finitely generated, then
Corollary~\ref{finid} also need not be true. Note that $\cI_A$ is
at most countable since idempotents are in one-to-one
correspondence with subgroups of $\bbZ^n$ generated by faces of
$S$ (the family of subgroups of $\bbZ^n$ is countable since each
of them can be generated by $n$ elements of $\bbZ^n$).
\begin{example}\rm
The light cone $r^2-p^2-q^2=0$ in $\bbR^3$ contains infinitely many
straight lines passing through points in $\bbZ^3$.
Let the semigroup $S\subset\bbZ^3$ be defined by inequalities
$r\geq0$, $r^2-p^2-q^2\geq0$. Then $\cI_A$ is infinite.
\end{example}
\noindent Perhaps, the most significant difference between
finitely and infinitely generated algebras occurs in the following
example.
\begin{example}\label{baapp}
\rm
Let $\al\in\bbR$ be irrational and put $S=\{(p,q)\in\bbZ^2:\,p+q\al>0\}$.
Then $A$ can be realized as an algebra of almost periodic
function which are bounded and analytic in the upper half-plane,
with the set of almost periods in the group $\{p-q\al:\,(p,q)\in\bbZ^2\}$;
$A$ is antisymmetric, has no orthogonal real measures, and is
a maximal subalgebra of $C(\bbT^2)$.
In \cite[Ch. 7]{Gam}, these algebras are studied in details.
\end{example}
\begin{proposition}\label{liesw}
A compact Lie group $G$ has the property {\rm(\sf{SW})} if and
only if its centre is finite.
\end{proposition}
\begin{proof}
If $Z(G)$ is not finite, then it contains a circle subgroup
$H\cong\bbT$; the family of all function in $C(G)$ that admit
analytic extension to $\bbD$ from each coset of $H$ is an
invariant algebra on $G$ which is not self-adjoint. The converse
follows from Proposition~\ref{cgace} and Corollary~\ref{cente}.
\end{proof}
A connected compact Lie group $G$ has a finite centre if and only
if $G$ is semisimple but in general its identity component may be
even abelian. Proposition~\ref{liesw} can be easily generalized to
all compact Hausdorff groups. A compact Hausdorff group is said to
be {\it profinite} if any neighbourhood of its identity contains a
normal subgroup of finite index.
\begin{theorem}\label{comsw}
A compact group $G$ has the property {\rm(\sf{SW})} if and only if
its centre is profinite.
\end{theorem}
\begin{proof}
Let $Z$ be the centre of $G$, $H$ be a closed normal subgroup of
$G$, and $\pi:\,G\to G/H=P$ be the canonical homomorphism. If $G$
satisfies ({\sf{SW}}), then $P$ also has this property (this is
evident). By the structure theorems, any neighbourhood of the
identity in $G$ includes a subgroup $H$ as above such that $P$ is
a Lie group. Since $\pi(Z)\subseteq Z_{P}(P)$, $\pi(Z)$ is finite
by Proposition~\ref{liesw}. Hence, $Z$ is profinite. Conversely,
let $G$ does not satisfy ({\sf{SW}}). By Lemma~\ref{redsw}, then
there exists a normal subgroup $H$ such that $P$ is a Lie group
that does not satisfy ({\sf{SW}}). Let $Z_P$ be the maximal torus
of $Z_{P}(P)$; by Proposition~\ref{liesw}, $Z_P$ is nontrivial. We
claim that $\pi(Z)\supseteq Z_P$. Indeed, if $G$ is a Lie group,
then for each $\xi\in\frz_P$, where $\frz_P$ is the Lie algebra of
$Z_P$, the affine space $\pi^{-1}(\xi)$ is $\Ad(G)$-invariant.
Hence, it contains an $\Ad(G)$-fixed point and the claim is true
for Lie groups. The general case reduces to this one since the
assumption $\pi(Z)\not\supseteq Z_P$ remains true if one replaces
$G$ by its quotient group over some sufficiently small normal
subgroup.
\end{proof}
Comparing Theorem~\ref{comsw} and the result of J. Wolf who gave
another characterization of these groups in \cite{Wo}, we get a
corollary: a compact Hausdorff group has a profinite centre if and
only if the image of each its one dimensional character is finite
(of course, this can be proved directly, starting with the case of
Lie groups).

A function algebra $A$ is called a {\it Dirichlet algebra} if
$A^\bot$ contains no nontrivial real orthogonal measure. In
\cite{Ri}, D. Rider proved that a compact group which admits an
antisymmetric Dirichlet algebra is connected and abelian. We prove
a more precise version of this theorem for Lie groups. First, we
describe a class of invariant algebras similar to the algebra of
example~\ref{baapp}. Let $v\in\bbR^n$ be such that
\begin{eqnarray}\label{dirri}
\scal{v}{x}\neq0\quad\text{for all}\quad x\in\bbZ^n\setminus\{0\}.
\end{eqnarray}
This happens if and only if components of $v$ are linear
independent over $\bbQ$.  Set
$$S_v=\{x\in\bbZ^n:\,\scal{v}{x}\geq0\}.$$
Let $A_v$ be the invariant algebra on $\bbT^n$ corresponding to
the semigroup $S_v$. Then $C^e=\bbR^+v$ and there are no other
nontrivial cones in $\frC_{S_v}$. If $\ga$ is the complex ray
corresponding to $v$, then $\ga(\bbC^+)$ is dense in $\cM_{A_v}$,
and $\ga(i\bbR)$ is an irrational winding of $\bbT^n$. Hence,
$A_v$ may be identified with the algebra $\{f\circ\ga:\,f\in
A_v\}$ of analytic functions on $\bbC^+$. Clearly, this algebra is
generated by exponents $e_x(z)=e^{iz\scal{v}{x}}$, where $x\in
S_v$, $z\in\bbC^+$. There are only two idempotents in $\cM_{A_v}$:
$e$ and $\eps$. Topologically, $\cM_{A_v}$ is the product
$\bbT^n\times[0,1]$ with identified points of the fibre
$\bbT^n\times\{0\}$.
\begin{theorem}
Let $G$ be a compact Lie group and $A$ be an antisymmetric
Dirichlet invariant algebra on $G$. Then $G=\bbT^n$ for some $n$.
The set $\cI_A$ is finite and linearly ordered; each group $G_j$
is a torus $\bbT^{n_j}$. Furthermore, if $j\leq k$ are consecutive
idempotents in $\cI_A$, then the algebra $A|_{kG_{j}}$ is closed
in $C(kG_j)$ and has the type $A_v$ described above, where the
vector $v\in\bbR^{n_j-n_k}$ satisfies {\rm(\ref{dirri})}.
\end{theorem}
\begin{proof}
Let $G_0$ be the identity component of $G$ and $L_{G_0}$ be the
averaging operator defined by (\ref{aveop}). If the algebra
$L_{G_0}A$ contains a nonconstant function, then it contains some
real function since $G/G_0$ is finite. Then $A$ cannot be
antisymmetric. If the averaging gives only constant functions,
then the difference between the Haar measures of $G$ and $G_0$ is
a real measure orthogonal to $A$. Hence, $A$ is not a Dirichlet
algebra if $G\neq G_0$.
By Theorem~\ref{antor} and Theorem~\ref{antis}, Haar measures on
$G$ and its maximal torus are representing for $\eps$, hence,
their difference is also a real orthogonal measure, which is
nontrivial if $G$ is not abelian. Therefore, $G$ is connected and
abelian. Since $G$ is a Lie group, $G=\bbT^n$ for some $n$.

Set $T^e=\frl^e+\Int C^e$. According to Theorem~\ref{anast}, there
exists a mapping $\eta:\,T\to\cM_A$ such that $f\circ\eta$ is
analytic and bounded on $T^e$ for all $f\in A$. If $\dim
\frl^e>1$, then for any point in $T^e$ there exist different
representing measures (for example, Poisson kernels on $\frl^e$
and on some line in $\frl^e$). For Dirichlet algebras, this is
impossible; hence, $\dim C^e=1$ and the closure of $S^e$ in
$\cM_A$ contains exactly one idempotent $j\neq e$. Therefore, each
chain of rays which starts at $e$ passes through $j$. Clearly, the
averaging of a Dirichlet invariant algebra over each closed
subgroup is again a Dirichlet one on the quotient group. Hence,
the arguments above can be applied to the algebra $R_jA$ on the
group $G/G_j$. Thus, Theorem~\ref{seque} implies that $\cI_A$ is
finite and linearly ordered.

Since $G_j$ is a p-set for $A$ by Proposition~\ref{idems}, 3), the
restriction of $A$ to $G_j$ is closed in $C(G_j)$. Obviously,
$A|_{G_j}$ is a Dirichlet invariant algebra. Hence, $G_j$ is a
torus.

Thus, it remains to prove that $A=A_v$, where $v\in\bbR^n$
satisfies {\rm(\ref{dirri})}, assuming that $\cI_A=\{e,\eps\}$.
Let $v$ generate $C^e$ and $\ga:\,\bbC^+\to\cM_A$ be the
corresponding complex ray. For any $z\in\bbC^+$, there exists the
unique representing measure for $\ga(z)$ concentrated on
$\ga(i\bbR)$ (the projection of the Poisson kernel on $\bbR$). The
support of any weakly limit point of these measures is contained
in the closure of $\ga(i\bbR)$. Since the Haar measure of $G$ is
the unique representing measure for $\eps$, $\ga(i\bbR)$ is dense
in $G$; i.e. it is an irrational winding of $G$. This property is
equivalent to {\rm(\ref{dirri})}. Set $S=\Sp A\subseteq\bbZ^n=\wh
G$. Then $S\subseteq S_v$ since $S^\star=C^e$. If
$\chi,-\chi\notin S$, then $\Re\chi\perp A$; hence,
\begin{eqnarray*}
S\mathbin\cup(-S)=\bbZ^n.
\end{eqnarray*}
Therefore, $S=S_v$.
\end{proof}
\section{Invariant algebras on spherical homogeneous spaces}
Let $M=K\backslash G$ be a right homogeneous space of a compact
connected Lie group $G$, where $K$ is the stable subgroup of the
base point $o$. In this section, we consider $G$-invariant
algebras on homogeneous spaces $M$, assuming that $M$ is {\it
multiplicity free}: the quasiregular representation of $G$ in
$C(M)$ contains every irreducible representation of $G$ with a
multiplicity $\leq1$. This remarkable class of homogeneous spaces
is characterized by each of the following properties (see, for
example, \cite{Vi}; if $G$ is not compact, then they are not
equivalent):
\begin{itemize}
\item[\rm 1)] $M$ is {\it commutative}, i.e., the algebra of all
invariant differential operators on $M$ is commutative; \item[\rm
2)] $(G,K)$ is a {\it Gelfand pair}, i.e., the convolution algebra
of all left and right $K$-invariant functions in $L^1(G)$ is
commutative; \item[\rm 3)] a generic $G$-orbit in the cotangent
bundle $T^*M$ is coisotropic; \item[\rm 4)] $M$ is weakly
commutative (i.e., the Poisson algebra $C(T^*M)$ is commutative);
\item[\rm 5)] $M$ is weakly symmetric; 
\item[\rm 6)] $M^\bbC=G^\bbC/K^\bbC$ is {\it spherical}, i.e., a
Borel subgroup of $G^\bbC$ has an open orbit in $M^\bbC$.
\end{itemize}
We shall say that $M$ is {\it spherical}. It is well known that
the group $G$ acting on itself by left and right translation is a
spherical homogeneous space, $G\times G/G$. In what follows, $M$
is supposed to be spherical and $A$ denotes a $G$-invariant closed
subalgebra of $C(M)$. Let $\wtd Z=Z_{\Diff(M)}(G)$ be the group of
all diffeomorphisms of $M$ that commute with every $g\in G$.
Clearly, any transformation in $\wtd Z$ keeps each irreducible
component of the quasiregular representation (which will be
denoted by $R$). The Schur lemma and the definition of the
multiplicity free spaces imply the following (well-known)
assertion.
\begin{lemma}\label{zinva}
If $M$ is spherical, then any $G$-invariant closed subspace of
$C(M)$ is $\wtd Z$-invariant. The group $\wtd Z$ is abelian. \qed
\end{lemma}
Transformations in $\wtd Z$ are induced by the left translations
$x\to gx$, where $x\in G$ and $g$ belongs to the normalizer of $K$
in $G$. Hence, $\wtd Z$ and $\wtd ZG=G\wtd Z$ are compact Lie
groups. Thus, studying $G$-invariant algebras on $M$, we may
assume without loss of generality that
\begin{eqnarray}\label{zsubg}
Z_{\Diff(M)}^0(G)\subseteq G,
\end{eqnarray}
where $Z_{\Diff(M)}^0(G)$ is the identity component of $\wtd Z$.
Let $Z=Z_G^0(G)$ be the identity component of the centre of $G$
and  $\si$ be the Haar measure of $Z$.
\begin{lemma}\label{pacom}
Let $A$ be antisymmetric, $H\subseteq G$ be a closed subgroup, and
$\nu$ be its Haar measure. If $H\supseteq Z$, then $R_\nu A=\bbC$,
where $\bbC$ denotes the set of constant functions on $M$.
\end{lemma}
\begin{proof}
It is sufficient to prove the lemma assuming $H=Z$: indeed, the
evident equality $R_\nu R_\si=R_\nu$ imply $R_\nu=R_\si$ since
both operators keep constant functions. By \cite[Theorem 8]{GL}
(whose most essential part was proved in \cite{La}), if $M$ admits
no one parameter group of transformations commuting with the
action of $G$, then each $G$-invariant algebra on $M$ is
self-adjoint with respect to the complex conjugation. Set $K'=KZ$;
it follows from (\ref{zsubg}) that the homogeneous space
$M'=K'\backslash G$ has this property. The space $R_\si A$
consists of all $Z$-invariant functions in $A$. Hence, it can be
considered as a subalgebra of $C(M')$. Since $A$ is antisymmetric,
we get $R_\si A=\bbC$.
\end{proof}
Let $H$ be a closed abelian subgroup and $\Sp_HA$ be the set of
weights of $R$ for $H$ in $A$. Since $A$ is an algebra, $\Sp_HA$
is a semigroup in $C(H)$. Hence, its closed linear span $B_H$ is
an invariant algebra on $H$.
\begin{lemma}\label{ancon}
Let $H$ and $B_H$ be as above, $\nu$ be the Haar measure of $H$.
Suppose that $H\supseteq Z$ and $A$ is antisymmetric.
Then $B_H$ is antisymmetric.
\end{lemma}
\begin{proof}
It is sufficient to prove {\rm(\ref{antab})}.
For $\chi\in\Sp_HA$, set
\begin{eqnarray}\label{pchi.}
P_\chi=\int R_z\ov{\chi(z)}\,d\nu(z).
\end{eqnarray}
If {\rm(\ref{antab})} is not true, then there exists a nontrivial
character $\chi\in\Sp_HA$ such that $\ov\chi\in\Sp_HA$. Then
\begin{eqnarray}\label{pchih}
P_\chi,P_{\ov\chi}\neq0\quad\mbox{in}\quad A.
\end{eqnarray}
Suppose $H=Z$. It follows from (\ref{pchih}) that there exist
nontrivial $G$-irreducible subspaces $U,V\subseteq A$ such that
\begin{eqnarray}\label{uvchi}
u(xz)=\chi(z)u(x),\quad  v(xz)=\ov{\chi(z)}v(x)
\end{eqnarray}
for all $u\in U$, $v\in V$, $z\in Z$, $x\in M$. The product $uv$
is $Z$-invariant. Hence $R_\si(uv)=uv$ and $uv=\const$ by
Lemma~\ref{pacom}. This implies $\dim (UV)\leq1$. Since $U$ and
$V$ are $G$-invariant and nontrivial, $UV\neq0$. Thus,
$\dim(UV)=1$. Since all functions in $U$ and $V$ are real
analytic, the assumption $u,v\neq0$, where $u\in U$ and $v\in V$,
implies $uv\neq0$. Consequently, $\dim U=\dim V=1$. It follows
that the semisimple part of $G$ acts on $U$ and $V$ trivially.
Thus, {\rm(\ref{uvchi})}, taken together with the structure
theorems and the additional (nonessential) condition
$u(o)=v(o)=1$, imply $v=\ov u$. Since $A$ is antisymmetric,
$u=v=\const$. Therefore, $\chi=\const$ contradictory to the
assumption.

Let $H\supset Z$. According to the proven above, if
$\chi\in\Sp_HA$ and $\ov\chi\in\Sp_HA$, then $\chi=1$ on $Z$.
Hence, any $\chi$-eigenfunction is $Z$-invariant. By
Lemma~\ref{pacom}, it is constant. Thus, $\chi$ is trivial.
\end{proof}
\begin{corollary}\label{mults}
Let $A$ be antisymmetric, $H$ and $\nu$ be as in Lemma~\ref{pacom}.
Then the mapping $R_\nu:\,A\to\bbC$ is a homomorphism.
\end{corollary}
\begin{proof}
By Lemma~\ref{ancon} and Theorem~\ref{antis}, $\nu$ is multiplicative
on the algebra $B_H$. Taken together with Lemma~\ref{pacom},
this proves the corollary.
\end{proof}
In fact, a stronger version of Corollary~\ref{mults} is true.
Let $H$ be as in Lemma~\ref{ancon} and $P$ be a closed subgroup of $H$
such that
\begin{eqnarray}\label{anmul}
P^\bot\mathbin\cap\Sp_HA=\{0\},
\end{eqnarray}
where $P^\bot=\{\chi\in\wh H:\,\chi\big|_P=1\}$ is the annihilator
of $P$ in $\wh H$, and let $\vk$ be the Haar measure of $P$.
\begin{corollary}\label{anhom}
If {\rm(\ref{anmul})} holds, then $R_\vk$ is a homomorphism of $A$
onto $\bbC$.
\end{corollary}
\begin{proof}
Due to Corollary~\ref{mults}, it is sufficient to prove that
$R_\vk=R_\nu$ in $A$, but this is a consequence of (\ref{anmul}).
\end{proof}
The following proposition generalizes  to the case of spherical
spaces an essential part of Theorem~\ref{antis}. We keep the
notation above. Let $\si$ be the Haar measure of $Z$ and $\mu$ be
the invariant measure on $M$ such that $\mu(M)=1$.
\begin{proposition}\label{sphan}
Following assertions are equivalent:
\begin{itemize}
\item[\rm 1)] $A$ is antisymmetric; \item[\rm 2)] $R_\si$ is a
homomorphism $A\to\bbC$; \item[\rm 3)] $\mu$ is multiplicative on
$A$; \item[\rm 4)] there is a $G$-fixed point in $\cM_A$.
\end{itemize}
\end{proposition}
\begin{proof}
The implication 1)$\ \Rightarrow\ $2) is a particular case of
Corollary~\ref{mults}. If $R_\si$ is a homomorphism $A\to\bbC$,
then for all $f\in A$ and $z\in Z$
\begin{eqnarray*}
\int f\,d\mu=\int R_zf\,d\mu=\int\int R_zf\,d\mu\,d\si(z)=\int
R_\si f\,d\mu= R_\si f(o),
\end{eqnarray*}
where the first equality holds since $\mu$ is invariant and the
last is true since $R_\si f=\const$. Hence, 2) implies 3). The
corresponding to $\mu$ point of $\cM_A$ is $G$-fixed. Conversely,
if $\nu\in\cM_A$ is $G$-fixed, then $\mu\in\cM_\nu$ since the set
$\cM_\nu$ is weakly compact, convex, and $G$-invariant. Thus, 3)$\
\Leftrightarrow\ $4). It remains to note that $M$ is a set of
antisymmetry if 3) is true since $M=\supp\mu$ and $\mu$ is a
representing measure.
\end{proof}
We use the assumption that $M$ is spherical only in the proof of
the implication 1)$\ \Rightarrow\ $2). Without this assumption,
the implication is false. For example, it is not true for adjoint
orbits $M=\Ad(G)v$ and algebras $A=P(M)$ (the closure of
polynomials on $M$), where $G=\SU(2)$, $v=i{\mathbf h}+r{\mathbf
e}$, $r>0$, and ${\mathbf h},{\mathbf e},{\mathbf f}$ is the
standard $\sll_2$-triple (cf. \cite{La} or \cite{GL}).

The canonical projection  $\pi:\,G\to M=K\backslash G$ induces an
embedding $\pi:\,C(M)\to C(G)$. In the following theorem, we
identify $A$ and $\pi(A)$, assume {\rm(\ref{zsubg})}, and keep the
notation of this section.
\begin{theorem}\label{acsph}
Let $M$ be spherical, $A$ be antisymmetric and $B$ be the closed
(bi)invariant algebra on $G$ generated by $A$. Then $B$ is
antisymmetric and the action of $G$ on $M$ extends to the action
of the semigroup $\cM_B$ on $\cM_A$ in such a way that for any
$\vf\in\cM_B$ the mapping $R_\vf$ defined by
\begin{eqnarray}\label{setrf}
R_\vf f(\psi)=f(\psi\vf)\quad\mbox{for all}\ \psi\in\cM_A,\ f\in A,
\end{eqnarray}
is an endomorphism of $A$. Moreover, $R_\vf=R_\mu$ for any
$\mu\in\cM_\vf$. If $\eps$ is zero of $\cM_B$, then $R_\eps
A=\bbC$.
\end{theorem}
\begin{proof}
Let $T$ be a maximal torus in $G$ and $B_T$ be as in
Lemma~\ref{ancon} (for $H=T$). Clearly, $\chi\in\Sp_TA$ if and
only if there exists $f\neq0$ in $A$ such that $f(mt)=f(m)\chi(t)$
for all $m\in M$, $t\in T$. It follows that $\Sp_TA$ is a
semigroup. As a closed linear space, $B$ is generated by products
of left shifts of functions in $A$. Since the left shifts commute
with the right ones, they do not change $\Sp_TA$; the products
also have spectrum in $\Sp_TA$ because $\Sp_TA$ is a semigroup.
Therefore,
\begin{eqnarray*}
\Sp_TB=\Sp_TA.
\end{eqnarray*}
Hence, $B_T$ coincides with the closure in $C(T)$ of the linear
span of characters in $\Sp_TA$. By Lemma~\ref{ancon} and
Theorem~\ref{antor}, $B$ is antisymmetric. It follows from
Lemma~\ref{chibo} that any $\vf\in\Hom(\Sp_TA,\bbD)$ defines an
endomorphism $R_\vf:\,A\to A$ by setting
\begin{eqnarray*}
R_\vf=\sum_{\chi\in\Sp_TA}\vf(\chi)P_\chi
\end{eqnarray*}
on $A_\fin$. Identifying $\Hom(\Sp_TA,\bbD)$ with $\cM_{B_T}$
according to Theorem~\ref{abelc}, we get an action of $\wh T$ on
$\cM_A$ that satisfies {\rm(\ref{setrf})}, where the right side
defines the left one.  This action extends to $\cM_B$ due to
Theorem~\ref{polar} and Theorem~\ref{Carta}, a). Clearly, for any
fixed $\psi\in\cM_B$ and all $f\in A$, the function $L_\psi
f(g)=f(\psi g)$, $g\in G$, belongs to $B$. Integrating by $g$ over
$\mu\in\cM_\vf$, we get $R_\vf=R_\mu$. The last assertion follows
from Theorem~\ref{antis}, c).
\end{proof}
Probably, the action of Theorem~\ref{acsph} is transitive (i.e.
$\cM_A=p\cM_B$ for any $p\in M$). Due to Theorem~\ref{surje}, this
is true if $A$ is the averaging of $B$ over $K$. For spheres
$\rU(n)/\rU(n-1)$, the latter was proved in \cite{Gi80} by methods
of harmonic analysis (which is rather simple on them). The
description of invariant algebras on these spheres can be found in
\cite{RuBS}, their maximal ideal spaces are characterized in
\cite{Ka}.

We conclude with an infinite dimensional version of the
Hilbert--Mumford criterion for commutative homogeneous spaces. It
admits several equivalent forms for finite dimensions. We
formulate a version which is a bit stronger than the classical one
and extends to the infinite dimensional case.

\smallskip
{\sl Let $G\subset\GL(V)$ be a compact connected Lie group, $\frg$
be its Lie algebra, and let $v\in V$. Suppose that
\begin{itemize}
\item[(N)] $p(v)=0$ for any $G$-invariant holomorphic
homogeneous polynomial $p$ of positive degree.
\end{itemize}
Then there exists $\xi\in\frg$ such that
\begin{eqnarray}\label{himum}
\lim_{t\to+\infty}e^{it\xi}v=0.
\end{eqnarray}
} In the standard statement for the field $\bbC$, the assertion
concerns an algebraic reductive group. Any such a group is the
complexification of a compact group $G$; if $T$ is a maximal torus
in $G$, then $G^\bbC$ has the Cartan decomposition $G^\bbC=GT^\bbC
G$. This makes it possible to reduce the problem to the compact
case. We omit the proof; the exposition below is self-contained.

\smallskip
Here is a simple example which shows that the criterion is not
longer true for Banach spaces, in the above form, in any sense. We
define a {\it homogeneous polynomial} $p$ of degree $n$ on a
Banach space $V$ by
\begin{eqnarray*}
p(v)=\phi(v,\dots,v),
\end{eqnarray*}
where $\phi$ is a continuous $n$-linear form on $V$. A {\it
polynomial} is a finite linear combination of homogeneous ones
(including constants), $\cP$ denotes the set of all polynomials on
$V$.
\begin{example}\label{nothm}
\rm
Let $A$ be any antisymmetric (bi)invariant algebra
on a compact Lie group $G$ and set
\begin{eqnarray*}
V=\left\{f\in A:\,f(\eps)=0\right\},
\end{eqnarray*}
where $\eps$ is the zero of $\cM_A$. Thus, $V$ is the unique
$G$-invariant maximal ideal in $A$. Clearly, $V$ is
$\frM_A$-invariant. For every $G$-invariant function $p$ on $V$,
any $v\in V$ and each ray $\ga\in\frR^e$ the function
$q(it)=p(R_{\ga(it)}v)$ does not depend on $t\in\bbR$ since
$\ga(it)$ lies in $G$. If $p$ is holomorphic, then $q(z)=p(v)$ for
all $z\in\bbC^+$. Therefore, $p(R_j v)=p(v)$ for $j=\ga(\infty)$.
The arguments above hold for the vector $R_jv$ and each ray
$\ga\in\frR^j$; applying them repeatedly, we get
\begin{eqnarray*}
p(v)=p(R_\eps v)=p(0)
\end{eqnarray*}
due to Theorem~\ref{seque}. It follows that each $G$-invariant
polynomial is constant; thus, (N) is true for all $v\in V$.

Let $\xi\in\frg$ satisfies {\rm(\ref{himum})} for a generic $v\in
V$ and $T$ be a torus which contains $\exp(\bbR\xi)$. Then
{\rm(\ref{himum})} holds for $P_\chi v$, there $\chi$ is a
character of $T$ and $P_\chi$ is the projection to the
$\chi$-isotypical component. Clearly, $\Sp_T v=\Sp_T A$ for a
generic $v\in A$; then {\rm(\ref{himum})} implies that $\xi\in
C^e$. Thus, $\exp(it\xi)$ defines  a ray $\ga\in\frR^e$ such that
$\lim_{t\to+\infty}\ga(t)=\eps$ but there exist algebras which
does not admit such a ray (see Example~\ref{notde}).\qed
\end{example}

Perhaps, this means that one parameter semigroup in a right
statement of any infinite dimensional version of the
Hilbert--Mumford criterion must be replaced by a chain of them.
We prove this for the spherical orbits (condition
{\rm(\ref{zsubg})} is not assumed). Operators $\xi\in\frg$ are
well-defined on the space $V_\fin$ of vectors $v$ such that the
linear span of $Gv$ is finite dimensional.
\begin{theorem}
Let $G$ be a compact connected Lie group acting strongly
continuously on a Banach space $V$. Suppose that $v\in V$
satisfies {\rm(N)} and that the orbit $M=Gv$ is a spherical
homogeneous space. Then there exist a non-increasing sequence
$V_1,\dots,V_n$ of closed subspaces of $V$, vectors $v_k\in V_k$,
where
$$k=1,\dots,n+1,\quad v_1=v,\ v_{n+1}=0,$$
and pairwise commuting $\xi_1,\,\dots,\xi_n\in\frg$ such that
$\xi_k$ and $e^{it\xi_k}$, $t\geq0$,  are well-defined on a dense
subspace of $V_k$ and
\begin{itemize}
\item[\rm1)] for any $t\geq0$, $e^{it\xi_k}$ extends to a
continuous operator $V_k\to V_k$; \item[\rm2)]
$\lim_{t\to+\infty}e^{it\xi_k}v_k=v_{k+1}$, $k=1,\dots,n$.
\end{itemize}
\end{theorem}
\begin{proof}
Let $A$ be the closure of $\cP$ in $C(M)$. Clearly, $A$ is an
invariant algebra on $M$. Averaging any $p\in\cP$ over $G$, we get
a polynomial (this is an easy consequence of the definition of a
homogeneous polynomial) which is necessarily $G$-invariant. It
follows from (N) that
\begin{eqnarray*}
&p(0)=\int_M p\,d\mu\quad\mbox{for all}\  p\in \cP,
\end{eqnarray*}
where $\mu$ is the normalized invariant measure on $M=Gv$ (note
that the converse is also true). Therefore, $\mu$ is
multiplicative on $\cP$, hence on $A$. Since $\supp\mu=M$, $A$ is
antisymmetric (see the proof of Theorem~\ref{antis}).

Let $K\subseteq G$ be the stable subgroup of $v$, $N$ be its
normalizer in $G$, and $H$ be a maximal torus in $N$. The group of
all transformations of $M$ which commute with $G$ may be
identified with $N/K$ acting on $M$ by right; by
Lemma~\ref{zinva}, $N/K$ is abelian. Let $Z$ be the identity
component of $N/K$ acting by right on $M$. Set $\td G=ZG$, $\td
H=ZH$ and let $\nu,\td\nu,\td\ka$ be the Haar measures of $H,\td
H,\td G$, respectively. Clearly, $\td H$ is a torus which contains
the identity component of the centre of $\td G$; hence we may
apply Lemma~\ref{ancon} and Corollary~\ref{mults} to $\td G$.
Since $\td Hv= Hv$, we get for any $p\in\cP$
\begin{eqnarray*}
\int_H p(hv)\,d\nu(h)=\int_{\td H} p(hv)\,d\td\nu(h)= \int_{\td G}
p(gv)\,d\td\ka(g)=\int_M p\,d\mu=p(0).
\end{eqnarray*}
Therefore, (N) holds for $H$ and $v$. Thus, we may assume $H=G$,
i.e., that $G$ is a torus. Set $V_0=V$, $v_1=v$,
\begin{eqnarray*}
S_1=\Sp_Hv_1=\{\chi\in\wh H:\,P_\chi v_1\neq0\},
\end{eqnarray*}
where $P_\chi$ is defined by {\rm(\ref{pchi.})}, and let $V_1$ be
the closure of $\sum_{\chi\in S_1}P_\chi V$. Then
\begin{eqnarray*}
\Sp_HV_1=\Sp_Hv_1\subseteq\Sp_HA.
\end{eqnarray*}
By Corollary~\ref{reext}, the representation of $G$ in $V_1$
extends to the representation of $\cM_A$ which is strongly
continuous. Let $\ga_1$ be a ray in $\cM_A$ such that
$\ga_1(0)v_1=v_1$ and $v_2=\ga_1(\infty)v_1\neq v_1$. Due to
Theorem~\ref{Liera}, there exists $\xi_1\in\frh$ such that
$\ga_1(t)v_1=\exp(it\xi_1)v_1$ for all $t\geq0$. Then (N) holds
for $v_2$ (see the beginning of Example~\ref{nothm}). Hence, the
invariant measure on $Hv_2$ is multiplicative on $\cP$ and on its
closure $A_2$. Applying this procedure repeatedly and replacing
$v_k$, $V_k$ by $v_{k+1}$, $V_{k+1}$, respectively, we get a
sequence of vectors $v_k$ such that $\dim Hv_k>\dim Hv_{k+1}$, if
$v_{k+1}\neq v_k$ (see Lemma~\ref{jkdim}). Let $v_n$ be the first
one that satisfies $v_{n+1}=v_{n+2}$. Since the closure $A_{n+1}$
of the restriction of $\cP$ to the torus $Hv_{n+1}$ is an
antisymmetric invariant algebra, Theorem~\ref{seque} and
Theorem~\ref{antis} imply $Hv_{n+1}=\{v_{n+1}\}$. If
$v_{n+1}\neq0$, then the averaging over $H$ of any continuous
linear functional that does not annihilate $v_{n+1}$ gives an
$H$-invariant linear functional with the same property,
contradictory to (N).
\end{proof}

\vbox{\vskip1cm
V.M. Gichev\\
gichev@ofim.oscsbras.ru\\
Omsk Branch of \\
Sobolev Institute of Mathematics\\
Pevtsova, 13, 644099\\
Omsk, Russia}
\end{document}